\documentclass[11pt]{article}

\usepackage{colordvi,epsfig,amsmath,amssymb}
\usepackage[table]{xcolor}   % for highlighting table cells
\usepackage[algosection,linesnumbered, ruled,noend]{algorithm2e}
\usepackage{hyperref}
\usepackage{verbatim,tikz}
\usepackage{caption,subcaption}   % Required for side-by-side figures with individual captions

\graphicspath{{figuras/}}  % Path to figures folder
\newcommand{\ignore}[1]{}
\newcommand{\R}{\mathbb{R}}
\newcommand{\N}{\mathbb{N}}
\newcommand{\K}{\mathcal{K}}
\newcommand{\qed}{\nobreak \ifvmode \relax \else
      \ifdim\lastskip<1.5em \hskip-\lastskip
      \hskip1.5em plus0em minus0.5em \fi \nobreak
      \vrule height0.75em width0.5em depth0.25em\fi}

\providecommand{\proj}[2]{\mbox{proj}_{#1}^{#2}}
\newcommand{\He}{\nabla^2 f} % D^2 f % \mathcal{B}
\newcommand{\norm}[1]{\lVert#1\rVert}
\newcommand{\dotprod}[1]{\left< #1\right>}
\newcommand{\Tr}[1]{\mbox{Tr}\left( #1\right)}
\newcommand{\quNac}{\mbox{\tt quNac }}

\def\proof{{\noindent{\bf Proof: }}}

\newtheorem{proposition}{Proposition}[section]
\newtheorem{lemma}{Lemma}[section]
\newtheorem{corollary}{Corollary}[section]

\title{Action constrained quasi-Newton methods} 
\author{R. M. Gower and J. Gondzio\footnote{School of Mathematics
and Maxwell Institute for Mathematical Sciences, 
The University of Edinburgh, corresponding author: gowerrobert@gmail.com}}

\begin{document}
\maketitle
%\tableofcontents

\begin{abstract}
At the heart of Newton based optimization methods
is a sequence of symmetric linear systems. Each consecutive system in this sequence is similar to the next, so solving them separately is a waste of computational effort. 
 Here we describe automatic preconditioning techniques for iterative methods for solving
 such sequences of systems by maintaining an estimate of the inverse system matrix. 
 We update the estimate of the inverse system matrix with quasi-Newton type formulas based on what we call an \emph{action constraint} instead of the secant equation. 
We implement the estimated inverses as preconditioners in a Newton-CG method and prove quadratic termination. Our implementation is the first parallel quasi-Newton preconditioners, in full and limited memory variants. Tests on logistic Support Vector Machine problems reveal that our method is very efficient, converging in wall clock time before a Newton-CG method without preconditioning. Further tests on a set of classic test problems reveal that the method is robust. The action constraint makes these updates flexible enough to mesh with trust-region and active set methods,  a flexibility that is not present in classic quasi-Newton methods.

%Here we describe a method that uses information made available from solving one system to update an estimate of the inverse matrix which is then used to precondition the next system. This information comes from matrix-free iterative solvers, where the system matrix is only applied as an operator. The operator's \emph{action} is sampled on multiple vectors while solving the system, and it is this sampled action that we use to update the estimated inverse system matrix.
%By modifying the quasi-Newton least-change formulation to incorporate action constraints, as opposed to the secant equation, we obtain a closed form solution of a family of updates.

\end{abstract}

\section{Introduction}
%%%%%%%%%%%%%%%%%%%%%%
\subsection{Motivation}
Second order methods for unconstrained nonlinear optimization display 
several advantages: they deliver a high accuracy of computations and 
enjoy a fast (quadratic) local convergence. However, these benefits may 
sometimes come at too high a cost. Indeed, evaluating the full Hessian 
and solving equations with it is sometimes very expensive and occasionally 
prohibitive. Several approaches have been designed over 
the years to remove some of the drawbacks of the second order methods 
while preserving their main advantages. Those include the {\it inexact} 
Newton methods~\cite{Dembo1982} and a family of {\it quasi-Newton} 
methods~\cite{Broyden1965, Broyden1970, Fletcher1970}.
  
%  
%  This information comes from matrix-free iterative solvers, where the system matrix is only applied as an operator. The operator's \emph{action} is sampled on multiple vectors while solving the system, and it is this sampled action that we use to update the estimated inverse system matrix.
%  
%  
  The inexact Newton method admits a (controlled) error 
in solving the Newton system and therefore allows to employ 
 matrix-free iterative solvers that only apply the system matrix as an operator.
These iterative methods only sample the action of the system operator, circumventing the cost of
calculating the entire Hessian matrix.
Quasi-Newton methods follow a completely different logic: they build 
an approximation of the inverse Hessian using low-rank updates derived 
from information on how the Hessian operates along a given direction.

The motivation behind this paper is to combine these two approaches: Use samples of
 the Hessian's action made available from an iterative solver to build an approximation to the inverse Hessian. This approximation is then used
to precondition and solve the subsequent Newton system, and the process is repeated.
%whose inexact solution produces samples of the subsequent Hessian's action which are used to refresh (and hopefully improve) the approximation of the inverse Hessian. 
The methods proposed in this paper and their analysis are based on the quasi-Newton literature.

The development of quasi-Newton methods was pioneered by Davidon
in the late 50's~\cite{Davidon1959} and culminated in the BFGS method,
named to honour the independent developments of Broyden~\cite{Broyden1965},
Fletcher~\cite{Fletcher1960}, Goldfarb~\cite{Goldfarb1970}
and Shanno\cite{Shanno1971} over the 60's and early 70's.
Nowadays, these methods are frequently referred to as members of the Broyden
family~\cite{Broyden1965, Broyden1970, Fletcher1970}.

Quasi-Newton methods obtain/improve an estimate $G_{k+1} \in S^n$ of the Hessian matrix 
$\He_{k+1} := \He (x_{k+1}),$ where $S^n$ is the set of symmetric matrices in $\R^{n \times n}$, $f\in C^2(\R^n)$ and $x_{k+1} \in \R^n.$
Their input is a previous estimate $G_k$ and a desired \emph{action} for the new estimate $G_{k+1}: \delta_k \rightarrow \gamma_k$, that is
\[ \gamma_k=  G_{k+1}\delta_k ,\]
where 
$\delta_k = x_{k+1} -x_k$ and $\gamma_k = \nabla f_{k+1}-\nabla f_k$. 
From the fundamental theorem of calculus
\[\gamma_k = \left(\int_{0}^1 \He(x_k + t\delta_k )dt\right) \delta_k, \]
so $G_{k+1}$ has the same action as $\int_{0}^1 \He(x_k + t\delta_k) dt$ applied to $\delta_k$. Alternatively, to obtain an estimate %$H_{k+1}\in S^n$
 of the (pseudo-)inverse Hessian, the action is inverted and imposed as $G_{k+1}: \gamma_k \rightarrow  \delta_k$. 

  This setup can produce approximate Hessians (or their inverse) from any observed action, in particular, when samples of the Hessian's action $d \rightarrow \He_{k+1}d$, with $d \in \R^n,$ are available. Though this limitation of incorporating only a 1-dimensional action is a hindrance when meshing quasi-Newton methods with inexact Newton methods because, in
contrast, inexact solvers make available the sampled action of the Hessian on a {\it subspace} (most often with more dimensions than one).
This mismatch has resulted in two strategies:
\begin{enumerate}
\item[(i)] Limit the inflow of new information to using only the action
of the Hessian on a {\it single} direction per iteration~\cite{Birgin2007,Inexact2006}. 
\item[(ii)] Use a basis for the subspace and associated Hessian's action, 
 to sequentially update the approximation~\cite{Morales2000a}.
This is costly and {\it cannot be parallelized.}
\end{enumerate}

We present a generalization of quasi-Newton methods which overcomes this drawback.

Instead of sampling the Hessian's action on a single direction, we sample it on a low dimensional subspace. This
guarantees a much faster influx of information and produces better approximations.
Using a set of directions at one time also allows us to
perform updates that exploit block-matrix operations which can
be executed in parallel.

Since the new methods exploit the Hessian's action along a set of directions,
we call them the {\it quasi-Newton Action Constrained} methods, \quNac
for short.

The motivation to develop \quNac comes from the need to solve large
and difficult problems. Therefore all computational aspects of the method  
are taken into serious consideration. In particular, we embed
\quNac into a Newton-CG scheme.
We discuss several variants
of a possible implementation of \quNac and provide preliminary
computational results which demonstrate its efficiency on non-trivial 
medium scale problems.

The next section contains the problem formulation and introduces 
the notation used in the paper. {From this initial motivation, we have broadened our scope to
 include preconditioning techniques for solving a sequence of (slowly) changing symmetric systems of equations
as opposed to focusing on a sequence of Newton systems.
Throughout the development we embrace two possible cases;
when \quNac approximations are developed either for estimating the system matrix or its inverse.

\subsection{Background}
%%%%%%%%%%%%%%%%%%%%

Consider the problem of sequentially solving in $d_k \in \R^n$ 
the symmetric systems
\begin{equation} \label{eq:seqQb}
   Q_k d_k = b_k, \quad \mbox{ for }k =1,2,\ldots,
\end{equation}
where $Q_k \in S^n$ and $b_k\in \R^n.$
Here the $Q_k$'s are ``slowly changing'' 
in the sense that $\norm{Q_{k+1} -Q_{k}}$ is relatively small in some 
matrix norm. We make no assumption on the $\{b_k\}$ sequence. 
Such slowly changing \emph{target matrices} $\{Q_k\}$ can arise from evaluating 
a continuous matrix field over neighboring points, such as is the case 
with the Hessian matrix in Newton type methods when step sizes are small. 
Sequences of symmetric systems also appear when solving nonlinear systems with the Newton-Raphson
method and the Jacobian is symmetric, such as discretizations
of the Nonlinear Schr\"odinger~\cite{sulem1999nonlinear} and the complex 
Ginzburg-Landau equation~\cite{ginzburg-landau}.

Solving a single system in~\eqref{eq:seqQb} through iterative methods involves calculating
$Q_{k+1} \mathcal{S}_k$, the action of $Q_{k+1}$ over a low dimensional 
\emph{sampling matrix} $\mathcal{S}_k \in \R^{n \times q}$, as opposed 
to requiring the entire matrix $Q_{k+1}$. This raises a question of how 
can one estimate the \emph{target matrix} $Q_{k+1}$, or its inverse, 
from this \emph{sampled action}. 

Our strategy is to maintain an \emph{estimate matrix} $G_k \in S^n$ of $Q_k$, 
and use the sampled action $\mathcal{S}_k \rightarrow Q_{k+1}\mathcal{S}_{k}$ 
to update $G_k$ and to produce a new estimate $G_{k+1} \in S^n$.
To determine a unique $G_{k+1}$, and exploit that $\norm{Q_{k+1}-Q_k}$ is small, 
we minimize $\norm{G_{k+1}-G_k}$ subject to an \emph{action} constraint
\[ G_{k+1} \mathcal{S}_k  = Q_{k+1}\mathcal{S}_k, \]
and a symmetry constraint
\[ G_{k+1} = G_{k+1}^T.\]
This is known as the \emph{least change} strategy in the quasi-Newton
literature, first proposed by Greenstadt in 1969~\cite{Greenstadt1969}. We henceforth refer to the problem of determining $G_{k+1}$ under these constraints
and the least change objective as the \emph{least change problem}.
As the constraint set $\{G \in \R^{n \times n} \, | \, G =G^T, \, G\mathcal{S}_k =Q_{k+1}\mathcal{S}_k\}$ is
a subspace of $\R^{n \times n}$, the resulting solution $G_{k+1}$  of the least change problem is a projection of $G_k$ onto this constraint set.
 This characterization as a projection is useful for including additional constraints as shown in the classic quasi-Newton setting by Dennis and Schnabel~\cite{Dennis1979}.
 
% We solve this problem using a weighted Frobenius norm in the objective 
% function and obtain a low rank update applied to $G_k$, akin to the
% quasi-Newton updates. 

The sampled action also offers information on the (pseudo-)inverse 
of $Q_{k+1}$ when it exists as 
\[Q_{k+1}^{-1}(Q_{k+1} \mathcal{S}_k) = \mathcal{S}_k.\]
Thus with an estimate $H_k\in \R^{n \times n}$ of the (pseudo-)inverse of $Q_{k}$,  %$H_k \in S^n$
 a new estimate can be obtained by minimizing the least change objective,
imposing the following action constraint
\[H_{k+1}(Q_{k+1} \mathcal{S}_k) = \mathcal{S}_k, \]
and the symmetry constraint. We use the same technique to calculate 
the direct or inverse estimate, the difference being which action we impose, 
$Q_{k+1} \mathcal{S}_k \rightarrow \mathcal{S}_k$ 
or $ \mathcal{S}_k \rightarrow Q_{k+1} \mathcal{S}_k.$

As our main application, we build estimates of inverse Hessian matrices 
to act as preconditioners in the Newton-CG method. In the unconstrained 
minimization of a function $f \in C^2(\R^n, \R)$, given an initial 
$x_0 \in \R^n,$ the Newton-CG method approximately solves systems,
\[\He_{k}d_k = - \nabla f_k,\]
using the Conjugate Gradient method~\cite{Hestenes1952}, where 
$\He_k:= \nabla ^2 f(x_k)$ is the Hessian matrix and $\nabla f_k := \nabla f(x_k)$, 
the gradient evaluated at $x_k \in \R^n.$ A line search is then performed 
to calculate a step size $a_k \in \R_+$ and iterate
\[x_{k+1} =x_k + a_k d_k.\]
In the Conjugate Gradient method, the action of the Hessian matrix 
is sampled on a low dimensional Krylov subspace. 
With this sampled  action we construct an estimate $G_k$ that is used to precondition 
the next Newton system $H_k \He_{k+1}d_{k+1} = -H_k\nabla f_{k+1}$.

%The remainder of the paper is as follows. After a description of previous work and connections to ours, we more on 
\subsection{Format of the paper}
After examining previous work and connections to our own in Section~\ref{sec:previouswork}, 
in Section~\ref{sec:COSYMOP} we solve the least change problem with a weighted Frobenius norm. 
Then we explore properties of the approximation matrices, such as sufficient conditions 
on the sampling matrix and target matrix to ensure the quadratic hereditary property and positive 
definiteness, both important in the context of preconditioning and in nonlinear optimization. This is followed by Proposition~\ref{prop:unravel} that shows when is the \quNac update equivalent to applying a sequence of rank-2 updates. This is used to establish the connection between sequential BFGS and DFP updates and \quNac updates.

We then specialize this updating scheme to Hessian 
matrices in Section~\ref{sec:invdirqunac} and develop a family of methods 
analogous to the Broyden family~\cite{Broyden1965}. In Section~\ref{sec:CG} 
we recap the Preconditioned Conjugate Gradients (PCG) method, followed  
by Section~\ref{sec:impinversequnac} where we detail a preconditioned 
Newton-CG method which employs \quNac in a full or limited memory variant that guarantee descent directions. We contrast our limited memory \quNac implementation to
Morales and Nocedal's L-BFGS preconditioner~\cite{Morales2000a}, showing that the former is a parallel version of the latter.
% 
%  Somewhat 
% surprisingly, we show how the limited memory implementation is a parallel 
% version of Morales and Nocedal's L-BFGS preconditioner~\cite{Morales2000a}.
The quadratic hereditary of this Newton-PCG method
is proved in Section~\ref{sec:qunacfull}, followed by promising numerical tests in Section~\ref{sec:tests}, comparing the new
method to Newton-CG, BFGS and L-BFGS on academic problems and regularized  logistic regression problems with real data.
 Finally we summarize our findings in Section~\ref{sec:conclusion}.

% In the next section, we solve a matrix optimization problem to obtain
% the $G_{k+1}$ estimate.

\subsection{Prior work and Connections}\label{sec:previouswork}

%This work draws upon selected results across numerical mathematics, extending a few, but more %importantly, uniting under a single perspective.

A member of the \quNac methods apparently first appeared in domain decomposition methods for solving PDE's~\cite{Mandel1993} where it is referred to as a balancing preconditioner. The domain decomposition methods give rise to a single large linear system which is block structured. After solving systems defined by the individual blocks, often in the least-squares sense, the balancing preconditioner aggregates these solutions into a symmetric preconditioner for the original large system. Our results enrich the balancing preconditioners by showing that they are a projection of a first guess preconditioner (the Neumman-Neumman preconditioner) onto the space of matrices with desirable properties (symmetric and having the same action as the (pseudo-)inverse over the direct sum of the nullspaces of the block matrices). Furthermore, we show that the balancing preconditioner is but one of a family of preconditioners that have these properties.
 % This connection to Balance preconditioner is by no means obviuos and is a cause for further investigation.

The balancing preconditioner has been taken out of the PDE context and tested as a general purpose preconditioner for solving a single linear system and systems with changing right hand side by Gratton, Sartenaer and Ilunga
~\cite{Gratton2011}. Gratton \emph{et al.} prove favourable spectral properties of the balancing preconditioner and study its relation to multiple BFGS updates. Our analysis of the quadratic hereditary property indicates how one might sequentially update a preconditioner
using the balancing preconditioner formula, which in turn allows us to extend the method to solving sequences of linear systems where the system matrix also changes.

%Gratton \emph{et al.} and 
% This work, rather, focuses on updating preconditioners while solving sequences of systems,
% where such results such quadratic hereditary are important. 

 The problem of solving sequences of linear systems has also been  addressed
by recycling Krylov subspace methods~\cite{Parks2006,Gaul2012,Loghin2006}
and in~\cite{Giraud2000Inc} when only the right-hand side changes.
In these methods, a selected Krylov subspace is retained from a previous
system solve that serves as an approximate eigenspace to improve
the conditioning of the next system.

Alternatively, updating a factorized preconditioner is possible,
such as partial LU decomposition for nonsymmetric systems~\cite{Tebbens2007}
and constraint preconditioners~\cite{Bellavia2014}.
  
Building a preconditioner through Frobenius norm probing~\cite{Huckle2007}
for a single linear system has a similar flavour to our preconditioning
method, where $H_{k+1}$ is obtained by approximately minimizing
$\norm{H_{k+1}Q_{k+1} -I}_F$ subject to an additional action constraint
that is incorporated into the objective function as a penalty. These aforementioned approaches,
and addressed problems, are notably distinct from ours. Rather, our setup
is heavily borrowed from that of quasi-Newton methods.

Schnabel~\cite{Schnabel1983} shows how to build estimate matrices that satisfy multiple secant equations, and in doing so, obtains generalizations of the Powell-Symmetric-Broyden (\emph{PSB}), BFGS and DFP formulas. He then goes on to show that
these generalizations are the solutions of the least change problem with a particular weighted Frobenius norm. By swapping multiple secant updates for an action constraint, Schnabel's generalized BFGS and DFP are equivalent to our inverse and direct \quNac method presented in Section~\ref{sec:invdirqunac}.

The least change problem was first proposed and solved for the standard quasi-Newton updates
\cite{Greenstadt1969,Goldfarb1970} but to the best of our knowledge this
paper is the first that solves the problem with a general action constraint and for any positive definite weighting matrix in the Frobenius norm.

Outside of the preconditioning literature, our proposed matrix optimization
problem has connections to low rank matrix completion~\cite{Candes2009}.
With a previous estimate $G_k =0$, one can view the action constraint
as a sampling of the target matrix through its action on a subspace.
The least change solution $G_{k+1}$ then leads to low rank solutions
of at most three times the number of columns in the sampling matrix.

\section{The quasi-Newton action constrained methods}
\subsection{The least change problem}
\label{sec:COSYMOP}
%%%%%%%%%%%%%%%%%%%

We now deduce the solution to the least change problem for a general action constraint and weighted Frobenius norm.
This includes and extends Schnabel's generalized BFGS, DFP and PSB methods~\cite{Schnabel1983}. 

Given an estimate matrix $G_k \in S^n$, our objective is to calculate an update matrix $E_k \in
S^n$ such that $G_k+E_k$ is an estimate of the target matrix $Q_{k+1} \in S^n.$
To ensure that the update matrix is the least change to $G_k$, it is obtained by minimizing a weighted Frobenius norm
\begin{equation} \label{eq:objfun}\norm{\mathcal{W}_k^{-1/2}E_k\mathcal{W}_k^{-1/2}}_{F}^2 :=
\Tr{\mathcal{W}_k^{-1}E_k\mathcal{W}_k^{-1}E_k^T},\end{equation}
where $\mathcal{W}_k \in S^n$ is a positive definite \emph{weighting} matrix. To impose that $G_{k+1}$ remains
symmetric, we
use a symmetry constraint
\begin{equation} \label{eq:symconst_1} E_k = E_k^T.\end{equation}
The action constraint is imposed as
 \begin{equation} \label{eq:actionconstgeniso} E_k \mathcal{S}_k = (Q_{k+1} -G_k)\mathcal{S}_k,\end{equation}
where $\mathcal{S}_k \in \R^{n \times q}$, $q$ an integer considerably smaller than $n$ and $\mathcal{S}_k$ is full rank.

Dropping the iteration index $k$, collecting the objective function~\eqref{eq:objfun}, symmetry
constraint~\eqref{eq:symconst_1} and the action constraint~\eqref{eq:actionconstgeniso} we have the \emph{least change problem}
that characterizes our update
\begin{align}
  \min_{E} & \, \frac{1}{2}\Tr{\mathcal{W}^{-1}E\mathcal{W}^{-1}E^T} \label{eq:obj}\\
 & \quad  E \mathcal{S} = R \mathcal{S} \label{eq:actionconstgenR}\\
& \quad E = E^T,  \label{eq:symconst}
\end{align}
where $R \in S^{n}$ is a given symmetric matrix.
 We now deduce the solution to the least change problem which is one of the central results of this article.
A key definition we repeatedly
use is
\[\proj{\mathcal{S}}{\mathcal{W}} :=  \mathcal{S}(\mathcal{S}^T\mathcal{W}\mathcal{S})^{-1}\mathcal{S}^T,\] thus
$\proj{\mathcal{S}}{\mathcal{W}}\mathcal{W}$ is an oblique projection onto
the space spanned by the columns of
 $\mathcal{S}.$ The following
demonstration is not necessary for the development of the remainder of the article, and the reader may jump ahead to the
solution~\eqref{eq:genupdate}.

The objective function of the least change problem~\eqref{eq:obj} is a convex quadratic function of $E$ and the constraints are
linear. Thus the solution is unique and characterized by the KKT conditions.
The Lagrangian of our least change problem is given by
\[\Phi(E,\Lambda,\Gamma) = \frac{1}{2}\Tr{\mathcal{W}^{-1}E\mathcal{W}^{-1}E^T} +\Tr{\Lambda^T(E-R)\mathcal{S}} + \Tr{\Gamma(E-E^T)},\]
where $\Lambda \in \R^{n \times q}$ and $\Gamma \in \R^{n\times n}.$
Differentiating (for a comprehensive list of formulas on matrix differentiation please
consult~\cite{Petersen2012}) in $E$ we have
\[D_E \Phi(E,\Lambda,\Gamma)= \mathcal{W}^{-1}E\mathcal{W}^{-1} + \Lambda \mathcal{S}^T +\Gamma^T -\Gamma. \]
Setting $D_E \Phi(E,\Lambda,\Gamma)$ to zero and isolating $E$ gives
\begin{equation} \label{eq:dem1}E =  \mathcal{W}(\Gamma - \Lambda \mathcal{S}^T -\Gamma^T)\mathcal{W}.\end{equation}
Using the symmetry constraint~\eqref{eq:symconst} of $E$ we find that
\[\Gamma -\Gamma^T =\frac{1}{2}\left(\Lambda \mathcal{S}^T - \mathcal{S} \Lambda^T \right). \]
Substituting back into~\eqref{eq:dem1} gives
\begin{equation} \label{eq:ELambda} E =  -\frac{1}{2}\mathcal{W}\left( \mathcal{S} \Lambda^T + \Lambda \mathcal{S}^T
\right)\mathcal{W}.
\end{equation}
The solution $E$ is now solely determined by $ \Lambda \mathcal{S}^T$, and we focus on obtaining this matrix.
Right multiplying by $\mathcal{S}$ and using the action constraint~\eqref{eq:actionconstgenR} then left multiplying by
$\mathcal{W}^{-1}$ gives
\begin{equation} \label{eq:rankrestrictL} \mathcal{W}^{-1}R\mathcal{S} =  -\frac{1}{2}\left( \mathcal{S} \Lambda^T + \Lambda
\mathcal{S}^T \right)\mathcal{W} \mathcal{S}.
\end{equation}
If the columns of $\mathcal{S}$ are linearly independent then $ \mathcal{S}^T \mathcal{W}\mathcal{S}$ is invertible.
Isolating $\Lambda$
\begin{equation}\label{eq:Lambdaiso}
\Lambda =  -\left( \mathcal{S} \Lambda^T \mathcal{W}\mathcal{S} +2\mathcal{W}^{-1}R\mathcal{S}\right) (\mathcal{S}^T
\mathcal{W} \mathcal{S})^{-1}.
\end{equation}
Right multiplying by $\mathcal{S}^T$ we find that
\begin{equation} \label{eq:LDproj} \Lambda \mathcal{S}^T = -\left( \mathcal{S} \Lambda^T \mathcal{W} +2\mathcal{W}^{-1}R
\right)
\proj{\mathcal{S}} {\mathcal{W}}.
\end{equation}
From~\eqref{eq:LDproj} we see that $\Lambda \mathcal{S}^T$ is equal to an unknown matrix times the
matrix $\proj{\mathcal{S}} {\mathcal{W}}$. This is a fact we shall use later on in the demonstration.
Left multiplying by $\mathcal{S}^T \mathcal{W}$ in~\eqref{eq:Lambdaiso}, we get
\[ \mathcal{S}^T \mathcal{W} \Lambda  = -\mathcal{S}^T \mathcal{W} \left( \mathcal{S} \Lambda^T \mathcal{W}\mathcal{S}
+2\mathcal{W}^{-1}R\mathcal{S}
\right)(\mathcal{S}^T \mathcal{W}\mathcal{S})^{-1} , \]
transposing
\[\Lambda^T \mathcal{W}\mathcal{S} =  -(\mathcal{S}^T \mathcal{W}\mathcal{S})^{-1}  \left( \mathcal{S}^T
\mathcal{W}\Lambda
\mathcal{S}^T +2\mathcal{S}^TR \mathcal{W}^{-1} \right) \mathcal{W}\mathcal{S}.\]
Substituting this into~\eqref{eq:Lambdaiso} we get
\begin{align*}
\Lambda & = \left(\phantom{\proj{\mathcal{S}} {\mathcal{W}}}\hspace{-1.0cm}
\mathcal{S}(\mathcal{S}^T
\mathcal{W}\mathcal{S})^{-1}  \left( \mathcal{S}^T \mathcal{W}\Lambda \mathcal{S}^T +2\mathcal{S}^TR \mathcal{W}^{-1}\right) \mathcal{W}\mathcal{S} - 2\mathcal{W}^{-1}R\mathcal{S} \right)(\mathcal{S}^T
\mathcal{W}\mathcal{S})^{-1}\\
& = \proj{\mathcal{S}} {\mathcal{W}} \mathcal{W}\Lambda +2\left(\proj{\mathcal{S}} {\mathcal{W}}R \mathcal{S} - \mathcal{W}^{-1}R\mathcal{S}\right) (\mathcal{S}^T \mathcal{W}\mathcal{S})^{-1}.
\end{align*}
 Right multiplying by $\mathcal{S}^T$ and isolating $\Lambda \mathcal{S}^T$ gives
\begin{align*}
 ( I -\proj{\mathcal{S}} {\mathcal{W}}\mathcal{W} ) \Lambda \mathcal{S}^T  &=2\left(\mathcal{S}(\mathcal{S}^T
\mathcal{W}
\mathcal{S})^{-1}\mathcal{S}^TR  - \mathcal{W}^{-1}R\right) \mathcal{S}(\mathcal{S}^T \mathcal{W} \mathcal{S})^{-1}
\mathcal{S}^T  \\
&=-2\left(I - \proj{\mathcal{S}} {\mathcal{W}}\mathcal{W}\right) \mathcal{W}^{-1}R  \proj{\mathcal{S}} {\mathcal{W}}.
\end{align*}
The above gives the $\left(I - \proj{\mathcal{S}} {\mathcal{W}}\mathcal{W}\right)$ projection of  $\Lambda
\mathcal{S}^T$.
It remains to find the $\proj{\mathcal{S}} {\mathcal{W}}\mathcal{W}$ projection of $\Lambda \mathcal{S}^T$. Decomposing
$\Lambda \mathcal{S}^T$ according to these projections we find
\begin{equation} \label{eq:Lformat}
 \Lambda \mathcal{S}^T = -2\left(I -  \proj{\mathcal{S}} {\mathcal{W}}\mathcal{W}\right) \mathcal{W}^{-1}R
\proj{\mathcal{S}} {\mathcal{W}} +
\proj{\mathcal{S}} {\mathcal{W}}\mathcal{W} \Lambda \mathcal{S}^T.
\end{equation}
  From~\eqref{eq:LDproj} we know that there exists $\Psi\in \R^{n\times n}$  such that $\Lambda \mathcal{S}^T = \Psi \proj{\mathcal{S}} {\mathcal{W}},$ thus
\begin{equation} \label{eq:Lformatpsi}
 \Lambda \mathcal{S}^T = -2\left(I -  \proj{\mathcal{S}} {\mathcal{W}}\mathcal{W}\right) \mathcal{W}^{-1}R
\proj{\mathcal{S}} {\mathcal{W}} +
\proj{\mathcal{S}} {\mathcal{W}}\mathcal{W}\Psi \proj{\mathcal{S}} {\mathcal{W}}.
\end{equation}
Inserting~\eqref{eq:Lformatpsi} into~\eqref{eq:LDproj}, after some elimination, we find that
\[ 2 \proj{\mathcal{S}} {\mathcal{W}}R \proj{\mathcal{S}} {\mathcal{W}} = -\proj{\mathcal{S}}
{\mathcal{W}}(\mathcal{W}\Psi+(\mathcal{W}\Psi)^T)\proj{\mathcal{S}} {\mathcal{W}}. \]
The solution is $\Psi = -\mathcal{W}^{-1}R,$ upto additions in the nullspace of $\mathcal{S}.$ This reduces~\eqref{eq:Lformat} to 
\[\Lambda \mathcal{S}^T =\left( \proj{\mathcal{S}} {\mathcal{W}}\mathcal{W} -2I\right) \mathcal{W}^{-1}R  \proj{\mathcal{S}}
{\mathcal{W}}. \]
Inserting the above in~\eqref{eq:ELambda} we obtain the solution
\begin{align} E &= -\frac{1}{2}\left(  \left( \mathcal{W}\proj{\mathcal{S}} {\mathcal{W}} -2I\right) R
\proj{\mathcal{S}}
{\mathcal{W}} \mathcal{W} +
\mathcal{W}\proj{\mathcal{S}} {\mathcal{W}}  R \left(\proj{\mathcal{S}} {\mathcal{W}}\mathcal{W} -2I\right) \right)
\nonumber\\
& = \mathcal{W}\proj{\mathcal{S}} {\mathcal{W}}  R\left(I -  \proj{\mathcal{S}} {\mathcal{W}} \mathcal{W}\right) + R
\proj{\mathcal{S}}
{\mathcal{W}} \mathcal{W}. \label{eq:genupdateE} \qed
\end{align}

 Picking up the iteration index $k$ again, identifying $R =Q_{k+1} -G_k,$ the projection of $G_{k}$ onto the subspace of symmetric matrices that satisfy the action constraint
is given by
\begin{equation}%\tag{\tt COSYsol}
\label{eq:genupdate}
\hspace{1.30cm} \boxed{ G_k+ E_k = Q_{k+1} +\left( I- \mathcal{W}_k\proj{\mathcal{S}_k}{\mathcal{W}_k}
\right)(G_k-Q_{k+1})\left(I -\proj{\mathcal{S}_k}{\mathcal{W}_k} \mathcal{W}_k\right)},% \hspace{0.0cm} \mbox{\tt
%COSYsol}
 \end{equation}
which is a rank-$3q$ update applied to $G_k$ that only requires knowing $Q_{k+1}\mathcal{S}_k$ and $\mathcal{W}_k \mathcal{S}_k.$ 
The updates~\eqref{eq:genupdate} include generalization of quasi-Newton methods, analogous to Schnabel's generalization with an action constraint in the place of multiple secant equations.
The generalized DFP and Powell-Symmetric-Broyden (\emph{PSB}) method are recovered by substituting $\mathcal{W}_k =Q_k$ and $\mathcal{W}_k =I$, respectively. The generalized BFGS method for estimating the inverse target matrix is recovered by substituting $\mathcal{W}_k=Q_k$
and swapping the occurrences of $Q_k\mathcal{S}_k$ and $\mathcal{S}_k$, so that $Q_k \mathcal{S}_k \rightarrow \mathcal{S}_k$ is the imposed action constraint.  Different from Schnabel's proof of the generalized BFGS updates, our solution does not assume that $G_k$ is invertible.

% \begin{equation}\tag{\tt COSYsol}\label{eq:genupdate}
%  E_k = \mathcal{W}_k\proj{\mathcal{S}_k}{\mathcal{W}_k} (Q_k -G_k)\left(I -  \proj{\mathcal{S}_k}{\mathcal{W}_k}
%\mathcal{W}_k\right) +
%(Q_k  -G_k) \proj{\mathcal{S}_k}{\mathcal{W}_k} \mathcal{W}_k,
%  \end{equation}

\begin{comment}
 Calculating the projections in~\eqref{eq:genupdate} requires the inversion of a $q \times q$ matrix, which is
computationally viable if $q$ is small. But this inversion is not necessary when the sampling basis is
$\mathcal{W}_k$-orthonormal, that is
 when $\mathcal{S}_k^T \mathcal{W}_k \mathcal{S}_k =I \in \R^{q \times q}.$ In this case~\eqref{eq:genupdate} results
in% an inexpensive update
\begin{equation}\label{eq:euclidupdate}E_k = Q_{k+1} - (I- \mathcal{S}_k\mathcal{S}_k^T) (G_k -Q_{k+1})\left(I -
\mathcal{S}_k\mathcal{S}_k^T\right). \end{equation}
\end{comment}

We now move on to sufficient conditions that guarantee the quadratic hereditary property and
positive definiteness of the resulting approximation matrix.

\subsection{The quadratic hereditary property}\label{sec:quadhered}

Iteratively updating an estimate $G_k$ using~\eqref{eq:genupdate}, we would like the estimate matrices to gradually converge to the target matrices. 
Though updating $G_k$ using~\eqref{eq:genupdate} results in an estimate with the desired action,  this update might have a destructive interference on the overall convergence. When the target matrices change little from one iteration to the next, the key to promoting convergence is guaranteeing that  the new estimate $G_{k+1}$ inherits the action of the previous estimate $G_k.$ 
 In the Proposition below, we prove that this convergence occurs if the target matrix is constant
for a number of iterations, say $\rho \in \N$ iterations.
% If the target matrix is the inverse Hessian matrix of a convex quadratic function, this means that %at iteration $\rho+1$ a full
%Newton step can be taken $d_{\rho}= -H_{\rho+1}\nabla f_{\rho+1}$
%and the solution is attained.

%This convergence is guaranteed by
%the \emph{quadratic Hereditary property}~\eqref{eq:hered}, which states that the estimate $G_{k+1}$
%satisfies the action constraint for all $\mathcal{S}_i$ of previous iterations $i=k,k-1 \ldots, 0.$

% This property helps promote convergence in unconstrained nonlinear optimization, in that, if the Hessian does
% not change much over a few iterations, the sequences of estimates can better approximate the sequence of Hessian
% matrices. In
% Proposition~\ref{prop:quadhered}, we prove a sufficient requirement for estimates updated via~\eqref{eq:genupdate} to
% satisfy the quadratic hereditary property.

For simplicity, assume that we have a sequence of full rank sampling matrices $\mathcal{S}_i \in \R^{n \times q_i}$ and $\rho,q_i \in \N$ for
$i=1,\ldots, \rho$ such that
$\sum_{i=1}^{\rho} q_i =n.$ 
\begin{proposition}[Quadratic Hereditary]\label{prop:quadhered}
 Let $G_0 \in S^n$ and $G_{k+1} = G_k +E_k$ defined by~\eqref{eq:genupdate} with $Q_k = Q \in S^n$ and
$\mathcal{W}_k \succ 0$ for $k=0, \ldots,
\rho.$ If $\mathcal{S}_k^T\mathcal{W}_{k}\mathcal{S}_i =0$ for every $i < k \leq \rho$ then
 \begin{equation} G_{k+1} \mathcal{S}_i =  Q\mathcal{S}_i , \quad
\mbox{for }i \leq k \leq \rho, \label{eq:hered}
\end{equation} %\tag{\texttt{Hered}$_{k+1}$}
and $G_{\rho+1} = Q.$
\end{proposition}
\proof The proof is by induction on $k$ that~\eqref{eq:hered} is true. 
 For $k=0$, our hypothesis becomes $
G_1 \mathcal{S}_0 = Q \mathcal{S}_0$ which is
equivalent to the action constraint~\eqref{eq:actionconstgeniso} with $k=0$. Suppose our hypothesis is true for $k-1$ and let us analyse the $k$ case.
For $i=k$, \eqref{eq:hered} is equivalent to the action constraint~\eqref{eq:actionconstgeniso}. For $i \leq k-1$, as $\mathcal{S}_k^T\mathcal{W}_{k}\mathcal{S}_i =0,$  we have
\[\proj{\mathcal{S}_{k}}{\mathcal{W}_{k}}\mathcal{W}_{k}\mathcal{S}_i =0.\]
Using~\eqref{eq:genupdate} to substitute $G_{k+1}$, we have
\begin{align*}
   G_{k+1}\mathcal{S}_i &= Q\mathcal{S}_i +\left( I- \mathcal{W}_k\proj{\mathcal{S}_k}{\mathcal{W}_k}
\right)(G_k-Q)\left(I -\proj{\mathcal{S}_k}{\mathcal{W}_k} \mathcal{W}_k\right)\mathcal{S}_i \\
&= Q\mathcal{S}_i +\left( I- \mathcal{W}_k\proj{\mathcal{S}_k}{\mathcal{W}_k}
\right)(G_k-Q)\mathcal{S}_i  \quad
\mbox{[by induction  $G_{k}\mathcal{S}_i = Q\mathcal{S}_i$, for $i\leq k$.]}  \\
&= Q\mathcal{S}_i.
\end{align*}
This concludes the induction. 

To prove $G_{\rho+1}=Q$, we need to show that the horizontal concatenation
 \[ \mathcal{S}_{1:\rho} : = \left[ \mathcal{S}_1 , \ldots,
\mathcal{S}_{\rho}\right], \]
is nonsingular.
To see this, let $\alpha_i \in \R^{q_i}$, for $i =0,\ldots, \rho$ be such that
\[\sum_{i=0}^{\rho}\mathcal{S}_i \alpha_i =0. \]
Left multiplying by $\alpha_{\rho}\mathcal{S}_{\rho}^T \mathcal{W}_{\rho}$ eliminates all terms except $\alpha_{\rho}^T\mathcal{S}_{\rho}^T \mathcal{W}_{\rho} \mathcal{S}_{\rho} \alpha_{\rho} =0,$ from which the positive definiteness of $\mathcal{W}_{\rho}$ and full rank of $\mathcal{S}_{\rho}$ implies
that $\alpha_{\rho}=0.$ The same procedure with $\alpha_{\rho-1}\mathcal{S}_{\rho-1}^T \mathcal{W}_{\rho-1}$ shows that $\alpha_{\rho-1} =0$
and so forth. Therefore, $\mathcal{S}_{1:\rho}$ has an inverse.
By induction~\eqref{eq:hered} is true for $k = \rho$, thus
\[   G_{\rho+1} \mathcal{S}_{1:\rho} =  Q \mathcal{S}_{1:\rho}.\]
Right multiplying the inverse of $\mathcal{S}_{1:\rho}$ on both sides shows that $G_{\rho+1} = Q. \qed$

To illustrate the proposition, consider the case where $W_k =I$ in~\eqref{eq:genupdate} which is a generalization of the PSB method. If the sampling
matrices $\mathcal{S}_i$ for $i=0, \ldots, k$ have mutually orthogonal columns, then 
Proposition~\ref{prop:quadhered} states that by updating using the PSB method the resulting $G_{k+1}$ satisfies the quadratic Hereditary property.
One way to achieve this would be to use residuals of a Krylov method to form the columns of the sampling matrices. 
Alternatively, if the weighting matrix satisfies the action constraint, then quadratic hereditary is guaranteed when the columns of the sampling matrix and resulting action matrix are orthogonal.

%\begin{corollary} \label{cor:quadhered}
% If $\mathcal{W}_{i}\mathcal{S}_i =Q\mathcal{S}_i$ for $i=0, \ldots,
%k$ and the columns of $[\mathcal{S}_0, \ldots, \mathcal{S}_{k}]$
%are $Q-$orthogonal %  , as will be the case in most of our applications,
%then due to Proposition~\ref{prop:quadhered}, the estimate matrix $G_{k+1}$ satisfies the quadratic Hereditary property.
%\end{corollary}

\begin{corollary} \label{cor:quadhered}
If  $\mathcal{S}_k^TQ\mathcal{S}_i =0$ for $i <k$ and $\mathcal{W}_{i}\mathcal{S}_i =Q\mathcal{S}_i$ for  $i\leq k$  then due to Proposition~\ref{prop:quadhered}, the estimate matrix $G_{k+1}$ satisfies the quadratic Hereditary property. \end{corollary}

The equivalent statements and proofs when the inverse action constraint $Q\mathcal{S}_k \rightarrow \mathcal{S}_k$ is imposed follow verbatim by swapping the labels of sampling matrix 
$\mathcal{S}_k$ and the sampled action $Q\mathcal{S}_k$. 
For example, after this label swap, Corollary~\ref{cor:quadhered} remains true though the weighting matrix need satisfy  $\mathcal{W}_{i}Q\mathcal{S}_i = \mathcal{S}_i$ and the resulting quadratic hereditary is $H_{k+1} Q\mathcal{S}_i =  \mathcal{S}_i$ for  $i\leq k$.

In the following section, we prove a sufficient condition for the solutions of the least change problem~\eqref{eq:genupdate} to be positive definite.

\subsection{Positive definiteness }\label{sec:posdef}
%{\bf dispense of indices?}

%Certain iterative solvers require positive definite preconditioners.
To apply the approximation matrix as a preconditioner, certain solvers require that it be positive definite.
Positive definiteness is also important in  unconstrained minimization: when
 we replace the Hessian matrix by an estimate matrix and solve the resulting quasi-Newton system, the search direction
is $d_k = -H_{k} \nabla f_k.$ If $H_k$ is positive definite and we are not at a stationary point $\nabla f_{k} \neq 0$ then $d_k$
is guaranteed to be a descent direction as
\[-\nabla f_k^{T}d_k = \nabla f_k^{T}H_k\nabla f_k > 0. \]
The next Lemma and Proposition are the main tools for proving positive definiteness of approximation matrices.

\begin{lemma}[Action Constrained Positive Definite Matrix] \label{lm:posdef} Let $P,A,B \in \R^{n \times n}$ where
 $A$ and $B$ are positive definite over Range$(P):=\{Px \, | \, x \in \R^n\}$ and Range$(I-P)$
respectively,
then the matrix
\[G = P^T A P + (I-P^T) B(I-P),\]
is positive definite.
\end{lemma}
\proof Let $x \in \R^n$, then
\begin{align*}
x^T Gx & = x^TP^TAP x+ x^T (I-P)^TB(I-P)x \geq 0.
\end{align*}
If $x^T Gx =0$ then $P x =0$ and $(I-P)x=0$ consequentially $x=Px + (I-P)x =0.$  $\qed$

% \begin{note} \label{nt:PAP}
%  When $P = \proj{\mathcal{S}}{A}A$, where $\mathcal{S} \in \R^{n \times q},$ then $P^T A P  = AP.$ 
% \end{note}

With Lemma~\ref{lm:posdef}, we characterize
 when a subset of estimate matrices that result
from~\eqref{eq:genupdate} are positive definite, namely those with a
weighting matrix that satisfies the action constraint
 $\mathcal{W}_k \mathcal{S}_k = Q_{k+1} \mathcal{S}_k.$ With such a weighting matrix, the
update~\eqref{eq:genupdate} takes the
form of the update~\eqref{eq:genquNac}, further down the page.
Such a weighting matrix always exists when $\mathcal{S}_k^TQ_{k+1} \mathcal{S}_k$ is positive
definite. To see this, let  $P = \proj{\mathcal{S}_k}{Q_{k+1}}Q_{k+1}$ and let
 \[\mathcal{W}_k =  Q_{k+1} P + (I-P)^T(I-P).\]
The projection matrix guarantees that $\mathcal{W}_k \mathcal{S}_k = Q_{k+1}\mathcal{S}_k$ and,
by noting that $Q_{k+1}P =P^T Q_{k+1} P,$  Lemma~\ref{lm:posdef} guarantees that the matrix $\mathcal{W}_k$  is positive definite.
 
% Note that only $Q_{k+1}\mathcal{S}_k$ is
%  required to compute the update and not the entire matrix $Q_{k+1}.$

\begin{proposition}[Positive Definite quNac] \label{prop:posdefk}
 If $G_0$ is positive definite and the product of the sampling matrix with the resulting action $\mathcal{S}_k^TQ_{k+1}\mathcal{S}_k $ is
positive definite  for $k =0, \ldots, \rho \in \N$ and
\begin{equation}\tag{\tt quNac} \label{eq:genquNac}
 G_{k+1} = Q_{k+1} \proj{\mathcal{S}_k}{Q_{k+1}}Q_{k+1} +\left(I- Q_{k+1} \proj{\mathcal{S}_k}{Q_{k+1}}\right)G_k
\left(I-
\proj{\mathcal{S}_k}{Q_{k+1}}Q_{k+1}\right),
\end{equation}
 then $G_{k}$ is positive definite for $k =0, \ldots, \rho+1$.
\end{proposition}

\proof By induction on $k$,  suppose that $G_k$ is positive definite. The first term on the right
hand side of~\eqref{eq:genquNac} can be
re-written as
\[ Q_{k+1} \proj{\mathcal{S}_k}{Q_{k+1}}Q_{k+1} = Q_{k+1}
\proj{\mathcal{S}_k}{Q_{k+1}}Q_{k+1}\proj{\mathcal{S}_k}{Q_{k+1}} Q_{k+1}.\]
In the context of Lemma~\ref{lm:posdef}, let $P = \proj{\mathcal{S}_k}{Q_{k+1}} Q_{k+1}$, $A =Q_{k+1}$ and $B=G_k$, and
by noting that Range$(P) = \mbox{Range}(\mathcal{S}_k)$ then
 $G_{k+1}$ is positive definite.$\qed$

 We call the estimates resulting from~\eqref{eq:genquNac} the \emph{quasi-Newton action constrained} estimates.
Different from~\eqref{eq:genupdate} which is a rank-$3q$ update, each
quNac estimate is a rank-$2q$ update. Next we prove an essential Lemma used to connect quNac methods to the BFGS and DFP methods.

From this point on, we apply~\eqref{eq:genquNac} as a function by
explicitly referring to the previous estimate and desired action 
$G_{k+1} =$\ref{eq:genquNac}$(G_k,\mathcal{S}_k \rightarrow Q_{k+1}\mathcal{S}_k)$. In particular,
in order the estimate an inverse matrix,  we apply the update $H_{k+1} =$\ref{eq:genquNac}$(H_k, Q_{k+1}\mathcal{S}_k \rightarrow \mathcal{S}_k)$ where the order of the action constraint has been switched.
Applying the positive definite Propositions to $H_{k+1}$ is simply a matter of switching the labels of $Q_{k+1}\mathcal{S}_k$ and $\mathcal{S}_k$ in the statements and proofs.
%The product of the sampling matrix with the resulting action is still $\mathcal{S}_k^TQ_{k+1}\mathcal{S}_k$, thus the Positive Definite proposition is applied verbatim.  

\subsection{Unravelling quNac into sequential rank 2 updates }
\label{sec:unravel}
Under orthogonality conditions between the columns of the sampling matrix and associated action, the rank-$2q$ ~\ref{eq:genquNac} update is equivalent to sequentially applying the quNac update built from the action on the $q$ individual columns of the sampling matrix. This has already been proved for the BFGS update in~\cite{Gratton2011}. We call this \emph{unravelling} the quNac update.

For this Proposition and henceforth, we say that $V,U \in \R^{n \times j}, j \in \N,$ are $A-$orthogonal, for $A\in S^n$, when $V^T A U = U^T A V = 0.$
\begin{proposition}[Unraveling] \label{prop:unravel} If the columns of $\mathcal{S}_k: =[s_1, \ldots,s_q ]$ are
$Q_{k+1}-$orthogonal, then \\
$G_{k+1}=\ref{eq:genquNac}(G_k,\mathcal{S}_k \rightarrow Q_{k+1}\mathcal{S}_k)$ is equal to $G_k^{q}$ where $G_k^1 := G_k$
and \[G_k^{i+1} = \ref{eq:genquNac}(G_k^i,s_i \rightarrow Q_{k+1}s_i), \quad \mbox{for }i =1,\ldots q.\]
\end{proposition}

\proof Borrowing Nocedal's notation~\cite{Nocedal1980} for
multiple BFGS updates, multiple quNac updates applied to $G_k$ to obtain $G_k^q$ is equivalent to
\begin{align}
 G_k^q &= (V_{1} \cdots V_{q})^T G_{k} (V_{1} \cdots V_{q}) \nonumber \\
&+ (V_{2} \cdots V_{q})^T Q_{k+1}\proj{s_1}{Q_{k+1}}Q_{k+1} (V_{2} \cdots V_{q}) \nonumber \\
&+ (V_{3} \cdots V_{q})^T Q_{k+1}\proj{s_2}{Q_{k+1}}Q_{k+1} (V_{3} \cdots V_{q}) \nonumber \\
&+ \cdots \nonumber \\
&+ Q_{k+1}\proj{s_q}{Q_{k+1}}Q_{k+1}, \label{eq:L-BFGS}
\end{align}
where $ V_i = I-\proj{s_i}{Q_{k+1}}Q_{k+1} $  for
$i=1,\ldots, q.$ As $s_i$ and $s_j$ are  $Q_{k+1}-$orthogonal for
$i\neq j$,
\begin{align*} V_i V_j &= (I- \proj{s_i}{Q_{k+1}}Q_{k+1} )(I- \proj{s_j}{Q_{k+1}}Q_{k+1} ) \\
 &= (I- \proj{s_j}{Q_{k+1}}Q_{k+1}  - \proj{s_i}{Q_{k+1}}Q_{k+1}  ) \\
&= (I- \proj{[s_j,s_i]}{Q_{k+1}}Q_{k+1}   ),
 \end{align*}
where $[s_j,s_i]$ is the column concatenation of $s_j$ and $s_i.$
This applied recursively yields
\begin{align*}
  (V_{i+1} \cdots &V_{q})^T Q_{k+1}\proj{s_i}{Q_{k+1}}Q_{k+1} (V_{i+1} \cdots V_{q})\\
& =  \left(I-Q_{k+1}\proj{[s_{i+1},\ldots, s_{q}]}{Q_{k+1}} \right) Q_{k+1}   \proj{s_{i}}{Q_{k+1}} Q_{k+1}\left(I-
\proj{[s_{i+1},\ldots, s_{q}]}{Q_{k+1}}Q_{k+1} \right) \\
&= Q_{k+1}\proj{s_i}{Q_{k+1}}Q_{k+1}.
\end{align*}
These observations applied to~\eqref{eq:L-BFGS} reveal
\begin{align*}
G_k^q &= Q_{k+1}\sum_{i=1}^q \proj{s_i}{Q_{k+1}}Q_{k+1} +\left(I- Q_{k+1} \proj{[s_{1},\ldots,
s_{q}]}{Q_{k+1}}\right)G_k \left(I-
\proj{[s_{1},\ldots, s_{q}]}{Q_{k+1}}Q_{k+1}\right)\\
&=  Q_{k+1} \proj{\mathcal{S}_k}{Q_{k+1}}Q_{k+1} +\left(I- Q_{k+1} \proj{\mathcal{S}_k}{Q_{k+1}}\right)G_k \left(I-
\proj{\mathcal{S}_k}{Q_{k+1}}Q_{k+1}\right)
\end{align*}
which is the quNac update $\ref{eq:genquNac}(G_k,\mathcal{S}_k \rightarrow Q_{k+1}\mathcal{S}_k). \qed$

Proposition~\ref{prop:unravel} is used to bridge quNac updates with sequentially applying Broyden family updates. Next we determine two practical quNac methods that generalize the DFP and BFGS methods.

\section{The \texttt{inverse} and \texttt{direct} quNac methods} \label{sec:invdirqunac}

Based on~\eqref{eq:genquNac}, we determine two methods for estimating the Hessian matrix $\nabla f_{k+1}$ and its (pseudo-)~inverse. The least change objective in the quNac framework can be
justified when $f$ is twice continuously differentiable, that is, $\He: x \rightarrow \He(x)$ is a
continuous matrix field.

With a given estimate $G_k \approx \He_k$, we define the \emph{direct} quNac update as  $G_{k+1} =
\ref{eq:genquNac}(G_k,\mathcal{S}_k~\rightarrow~\He_{k+1}~\mathcal{S}_k).$
Positive definiteness is guaranteed by Proposition~\ref{prop:posdefk} when $G_k \succ 0$ and
when $\mathcal{S}_k^T \He_{k+1}\mathcal{S}_k~\succ~0.$
 Using the Woodbury formula~\cite{Woodbury1950}, in the Appendix~\ref{sec:DFPinv} we show that much like the DFP method,
one can update the inverse when $H_k = G_k^{-1}$ exists and work solely with $H_k$ through the formula
\begin{equation}  %\tag{DirQ}
 \label{eq:directquNacH}
H_{k+1} = H_k+ \proj{\mathcal{S}_k}{\He_{k+1}}-  H_k\He_{k+1} \proj{\mathcal{S}_k}{\He_{k+1} H_k \He_{k+1}}\He_{k+1}
H_k.
\end{equation}

Alternatively, we can use the \ref{eq:genquNac} update to estimate the inverse Hessian without the need to go through the Woodbury formula. 
To build an estimate matrix $H_{k+1}\in S^n$ of the inverse Hessian with the appropriate action
 $H_{k+1}: \He_{k+1}\mathcal{S}_k \rightarrow \mathcal{S}_k$, we simply invert the order of the arguments
$\mathcal{S}_k$
and $\He_{k+1}\mathcal{S}_k$ in the \ref{eq:genquNac} function so that $H_{k+1} =
\ref{eq:genquNac}(H_k,\He_{k+1}\mathcal{S}_k \rightarrow \mathcal{S}_k ).$
This results in the \emph{inverse} quNac update
\begin{equation} %\tag{InvQ}
\label{eq:invquNac} H_{k+1} =
\proj{\mathcal{S}_k}{ \He_{k+1}} + \left(I-\proj{\mathcal{S}_k}{ \He_{k+1}}\He_{k+1}\right)H_k\left(I- \He_{k+1}
\proj{\mathcal{S}_k}{
\He_{k+1}}\right).
\end{equation}
In this inverse perspective, $\He_{k+1}\mathcal{S}_k$ is the sampling matrix and $\mathcal{S}_k$ the resulting action.
Positive definiteness of $H_{k+1}$ follows by Proposition~\ref{prop:posdefk} when $H_k \succ 0$ and when the product of the sampling matrix and associated action is positive definite, that is, when $\mathcal{S}_k^T\He_{k+1}\mathcal{S}_k \succ 0$.

The BFGS and DFP methods are instances of the inverse and direct quNac, respectively.
 When $\mathcal{S}_k = s \in \R^n$ is
comprised of a single column, then the inverse (direct) quNac update is equivalent to applying a BFGS (DFP) update with the action $\He_{k+1}s \rightarrow s$ $\left(s \rightarrow \He_{k+1}s\right)$ which can be re-written as
\begin{align*} \label{eq:BFGS}
H_{k+1} &= \frac{s s^T}{s^T \He_{k+1}s} +
\left(I-\frac{s s^T\He_{k+1}}{s^T \He_{k+1}s} \right)H_k\left(I -\frac{\He_{k+1}s s^T}{s^T \He_{k+1}s}\right)\\
& = \proj{s}{\He_{k+1}} + \left(I-\proj{s}{\He_{k+1}}\He_{k+1} \right)H_k\left(I -\He_{k+1}\proj{s}{\He_{k+1}}\right).
\end{align*}
That is, applying the BFGS and DFP update using the pair $\delta_k,\gamma_k \in \R^n$ is equivalent to applying the
update $ \ref{eq:genquNac}(H_k,\gamma_k \rightarrow \delta_k )$ and $ \ref{eq:genquNac}(G_k,\delta_k \rightarrow \gamma_k )$, respectively. Thus we can
apply Propositions~\ref{prop:posdefk} and~\ref{prop:quadhered} to show that the resulting estimate is positive
definite when $H_k \succ 0$, $\gamma_k^T \delta_k >0$ and quadratic Hereditary holds when $\{\delta_1, \ldots, \delta_k
\}$ are $Q-$orthogonal where $Q$ is the constant Hessian matrix. These sufficient conditions are well known for the BFGS and DFP
methods, but it is nice to see how they are derived using the same tools for quNac methods.

Furthermore, when the columns of $\mathcal{S}_k$ are $\He_{k+1}-$orthogonal,  then according to
Proposition~\ref{prop:unravel} applying the inverse (direct) quNac update  is equivalent to sequentially applying BFGS (DFP)
updates built from the $i$th column of $\mathcal{S}_k$ and $\He_{k+1}\mathcal{S}_k$, for $i =1,\ldots, q.$ We use this
observation to implement a new {\it parallelizable} method for applying a L-BFGS preconditioner.

 We now digress from the main flow of the article to show that,
much like the Broyden family, the \texttt{direct} and \texttt{inverse} quNac methods can be combined to generate a family of methods.

\subsection{A Family of quNac methods} \label{sec:family}
%{\bf Remove this subsection barrier and make into one big section? But then the user feels he must read the entire
%section, where according the hyper-contents, this subsection can be dispensed with.}
We can update a given $H_k$ estimate using  a combination
\begin{align*}
H^{\lambda}_{k+1} &= \lambda_k H_{k+1}^{D} +(1-\lambda_k)H_{k+1}^{I},
\end{align*}
where $H_{k+1}^I$ and $H_{k+1}^D$ are given by the inverse~\eqref{eq:invquNac} and direct~\eqref{eq:directquNacH} estimate,
respectively, and $\lambda_k \in [0 ,\, 1].$
Manipulating the formulas for $H_{k+1}^{D}$ and $H_{k+1}^{I}$ we find
\begin{align}
H^{\lambda}_{k+1}&= H_{k+1}^{I}  +\lambda_k\proj{\mathcal{S}_k}{\He_{k+1}}\He_{k+1} H_k\left(I -\He_{k+1}
\proj{\mathcal{S}_k}{ \He_{k+1}}
\right)
\\
&+\lambda H_k\He_{k+1}\left( \proj{\mathcal{S}_k}{ \He_{k+1}} -\proj{\mathcal{S}_k}{\He_{k+1} H_k \He_{k+1}}\He_{k+1}
H_k \right) \nonumber
\\
&= H_{k+1}^{I} -\lambda_k V_kV_k^T,\label{eq:quNaccombo}
\end{align}
where
\[V_k = \left(\proj{\mathcal{S}_k}{\He_{k+1}}\He_{k+1}  - I\right)H_k\He_{k+1}\mathcal{S}_k(\mathcal{S}_k^T\He_{k+1} H_k
\He_{k+1}\mathcal{S}_k)^{-1/2} \in \R^{n \times q},\]
thus analogously to the Broyden family, each member of the quNac family is at most a rank-$q$ matrix in distance from
each other.
When $H_{k+1}^{D}$ and $H_{k+1}^{I}$ are positive definite, then so is $H^{\lambda}_{k+1}$ as it is a positive sum of
two positive definite matrices.

The resulting $H^{\lambda}_{k+1}$ also satisfies the action constraint as
\begin{equation} \label{eq:convactionconst}  \left(\lambda_k H_k^{D} +(1-\lambda_k)H_k^{I}\right) \He_{k+1}
\mathcal{S}_k =
\lambda \mathcal{S}_k +(1-\lambda)\mathcal{S}_k
=\mathcal{S}_k.\end{equation}
When the quadratic Hereditary property holds for $H_{k+1}^{D}$ and $H_{k+1}^{I}$, it also holds for $H^{\lambda}_{k+1}$
using the same observation as in~\eqref{eq:convactionconst} though with $\mathcal{S}_i$ for $i =1,\ldots, k,$ in the place of $\mathcal{S}_k.$

To implementing a Newton-CG method with a \quNac preconditioner we need the details of the PCG method. 
 Readers familiar with the PCG method can jump to the Restarting Preconditioner Lemma~\ref{lemma:restart}.

\section{Conjugate Gradients} \label{sec:CG}
%Quick and painful introduction to Conjugate Gradients emphasizing the mathematical properties, i.e., all projections
%arev%named explicitly.

The conjugate gradients method, developed by Magnus Hestenes and Eduard Stiefel~\cite{Hestenes1952}, is an iterative method for finding the solution to
\begin{equation} \label{eq:phiquad}  \min_x \phi(x) := \min_x \frac{1}{2}x^TQx -x^Tb,\end{equation}
where $x,b\in \R^n$ and $Q \in S^n$ is a positive definite matrix which guarantees that the critical point defined by 
\begin{equation} \label{eq:residual}
\nabla \phi(x) = Qx - b =0,
\end{equation}
 is the unique solution. With a given $x_0 \in \R^n$,
the method iteratively finds $x_{k}$, the minimum of
$\phi(x)$ restricted to $x_0\oplus \K_k$, where $\K_k = \mbox{span}\left\{\nabla \phi(x_0), Q\nabla \phi(x_0), \ldots, Q^{k-1}\nabla \phi(x_0)\right\}$ is
the $k$th \emph{Krylov subspace}.  This
construction implies that if $v \in \K_k$ then $Qv \in \K_{k+1}.
$ The Krylov subspaces are nested, in that $\K_{k} \subset \K_{k+1}$, thus each $x_{k+1}$ tends to be an
improvement
over the previous $x_{k}.$
As $x_{k}$ is a constrained optima, the gradient $r_{k}:= \nabla \phi(x_k)$, which is the \emph{residual}  in
equation~\eqref{eq:residual} at $x_k$, is in $\K_{k}^{\perp}$, the
orthogonal complement of $\K_{k}.$

The CG method searches the Krylov spaces by using $Q-$orthogonal directions, which are also known as the conjugate directions.
The first conjugate direction is set to
$p_0:=-r_0$. An exact line search is then performed with $\alpha_0 := \arg \min\{\alpha \, | \, \phi(x_0 +\alpha
p_0)\}$
to obtain a new iterate $x_1 =x_0 + \alpha_0 p_0.$ For this reason $r_1$ is orthogonal to $\K_1 =
$span$\{p_0\}.$
Then recursively from
$x_k$, a conjugate direction in $\K_{k+1}$ is determined
by applying the Gram-Schmidt orthogonalization process with inner product $\dotprod{\cdot, \cdot}_Q$ to $-r_k$, 
\begin{align} p_k &= -r_k +\frac{\dotprod{r_k,p_{k-1}}_Q}{\dotprod{p_{k-1},p_{k-1}}_Q} p_{k-1}. \label{eq:dupdate}
 \end{align}
Only the component of $r_k$ in the $p_{k-1}$ direction is removed as $r_k \in \K_{k}^{\perp} \subset
(Q\K_{k-1})^{\perp}$ which guarantees that the inner product of $r_k$ with each  $Qp_0, \ldots, Qp_{k-2}$ is zero.
An exact line search over $p_k$ is then performed to find $x_{k+1}$
%  x update
\begin{align}
x_{k+1} &= x_k +\alpha_k p_k, \label{eq:xupdate}
\end{align}
where $\alpha_k = - \dotprod{ r_k, p_k}/\dotprod{ p_k,p_k}_Q.$
% Don't need this anymore -->
%Note that substituting $p_k$ for any scalar multiple of $p_k$ produces the same sequence of steps. We make use of this
%fact in the next section.
Finally, as $\phi(x)$ is a quadratic function,
the gradient can be calculated iteratively
\begin{align}
 r_{k+1} & = r_k +\alpha_k Q p_k. \label{eq:rupdate} 
%&=  r_k -\frac{\dotprod{ r_k, p_k}}{\dotprod{ p_k,p_k}_Q } Q p_k. 
\end{align}
If a preconditioner $M\in S^n$ with $M \succ 0$ is used, in other words, if an equivalent positive definite system to
$M^{-1}Qx =
M^{-1}b$ is solved, then the
Gram-Schmidt process is applied to the sequence $M^{-1}r_k$ instead of $r_k$ resulting in
\begin{align}
p_0 &= -M^{-1}r_0, \label{eq:Md0update}\\
 p_{k} & = -M^{-1}r_k +\frac{\dotprod{M^{-1}r_k,p_{k-1}}_Q}{\dotprod{p_{k-1},p_{k-1}}_Q} p_{k-1}, \quad k>0.
\label{eq:Mdupdate}
\end{align}

% \Rob{Not happy with this next paragraph and how it affects the flow. How can I justify appearance of the following proposition, without going into details of the quNac Newton-PCG method, which would be out of place.}
Before moving on, we need a Lemma that is fundamental in proving the quadratic Hereditary property of our forthcoming Newton-PCG implementation.
The Lemma establishes sufficient conditions on the preconditioner and a new starting point such that after stopping then starting the PCG method at this new point,
 the PCG method continues to build $Q$-orthogonal search directions.

% 
% In both cases, with or without a preconditioner, the conjugate search directions $p_1, \ldots, p_k$ form a  $Q$-orthogonal basis for a Krylov subspace.
% Iin our forthcoming Newthon-PCG method, the columns of $\mathcal{S}_k$ are  the conjugate search directions stored from applying the Newton-PCG method to
%  $\He_{k} d_k= -\nabla f_{k}$ with a preconditioner $H_k$ previously calculated. When proving the quadratic hereditary of the updated preconditioner $H_{k+1}$ later on, we need the following Lemma.

% Previous Explanation:
% In calculating quNac estimates, a few iterations of the PCG method will be applied to the Newton system $\He_{k+1} d= -\nabla f_{k+1}$ to obtain a Krylov subspace that
%  calculate the columns of $\mathcal{S}_k$ on the $(k+1)$th iteration.
% But as the PCG method is restarted on a new system at each iteration $k$,
%  we cannot guarantee that the columns of $\mathcal{S}_{k}$ will be $Q-$orthogonal to those
% of $\mathcal{S}_i$, for $i<k,$  as required by the quadratic Hereditary proposition.
% Though preconditioners that enjoy the quadratic Hereditary property guarantee $Q-$orthogonality across restarts, as shown in the next Lemma.

\begin{lemma}[Restarting Preconditioner] \label{lemma:restart}
 Let $p_{0} \ldots p_{k-1}$ be a set of $Q-$orthogonal directions. Let $\bar{x}_0 \in \R^n$ with gradient $\nabla
f(\bar{x}_0)$  such that $p_j^T \nabla f(\bar{x}_0) =0,$ for $j=1,\ldots, k-1.$ Let $M\in \R^n$ be a symmetric
positive
definite matrix such that
\begin{equation} \label{eq:dieigen} M^{-1} Q p_j = p_j, \quad \mbox{for } j =1,\ldots, k-1. \end{equation}
Then by executing $t$ iterations of the PCG method on the system $Qx =b$, where $k+t+1 \leq n,$ with initial point $\bar{x}_0$ and $M^{-1}$ as
a preconditioner, the conjugate directions calculated, namely
 $\bar{p}_0, \ldots, \bar{p}_t$,  are such that 
\[\{p_{0} \ldots p_{k-1}, \bar{p}_0, \ldots, \bar{p}_t\},\] is a $Q-$orthogonal set. 
\end{lemma}
\proof
Let $ \bar{r}_0 , \ldots, \bar{r}_t$ be the residual vectors associated with the conjugate directions $\bar{p}_0, \ldots, \bar{p}_t,$ where $\bar{r}_0 := \nabla f(\bar{x}_0).$
We use induction on $t$, where our induction hypothesis is that  $\bar{p}_i^T Qp_j  =0$ and $\bar{r}_i^T p_j=0$ for
$1\leq j \leq k-1$ and $0 \leq i \leq t$. For $t=0$, as $\bar{p}_0 = -M^{-1}\bar{r}_0$,
 \begin{align*}
  \bar{p}_{0}^T Q p_j   &=    -\bar{r}_{0}^T M^{-1}  Q p_j \quad \left(\mbox{using~\eqref{eq:dieigen}}\right) \\
 &= -\bar{r}_0^T p_j  = 0, \quad\mbox{for } j =1,\ldots, k-1.
 \end{align*}
 Supposing the induction hypothesis is true for $t-1$ and all $0\leq j\leq k-1$, using~\eqref{eq:rupdate} to calculate
the next residual $\bar{r}_{t}$, then by induction
%step using~\eqref{eq:dupdate} and~\eqref{eq:rupdate}
\begin{align*}
\bar{r}_{t}^T p_j  &= \bar{r}_{t-1}^T p_j -\frac{\dotprod{ \bar{r}_{t-1}, \bar{p}_{t-1}}}{\dotprod{
\bar{p}_{t-1},\bar{p}_{t-1}}_Q }  \bar{p}_{t-1}^T Q p_j\\
&=  \bar{r}_{t-1}^Tp_j \\
&= 0.
\end{align*}
Using~\eqref{eq:Mdupdate} to substitute $\bar{p}_{t}$ 
\begin{align*}
\bar{p}_{t}^T Qp_j  &= -  \bar{r}_t^T M^{-1} Q p_j
+\frac{\dotprod{M^{-1}\bar{r}_t,\bar{p}_{t-1}}_Q}{\dotprod{\bar{p}_{t-1},\bar{p}_{t-1}}_Q} \bar{p}_{t-1}^T Qp_j\\
&= -\bar{r}_{t} ^T   M^{-1} Qp_j \quad \left(\mbox{applying~\eqref{eq:dieigen}}\right) \\
&= -\bar{r}_{t}^T p_j \\
&= 0, \quad \mbox{ for } j=1,\ldots, k-1. \qed
\end{align*}
We refer to $\bar{x}_0$ and $M^{-1}$ of Lemma~\ref{lemma:restart} as a restart point and restarting preconditioner,
respectively.

For further reading on the Preconditioned Conjugate Gradients (\emph{PCG}) method see~\cite{Shewchuk1994} for a
pedagogic explanation and~\cite{Gower2014a} for a description that uses oblique projections.

\section{Implementing a Newton-PCG quNac method}\label{sec:impinversequnac}
% Implementing a parallel BFGS preconditioned based on the quNac update.

%{\bf Why is $H_k$ a favorable pre-conditioner and how in the quadratic case all previous PCG direction are e-vec with
%e-vals 1.}\\

We use the \texttt{inverse} quNac formula~\eqref{eq:invquNac} to update a preconditioner within a Newton-PCG method for  finding local minima of $f \in C^2(\R^n),$ where $f$ is possibly non-convex,  see Algorithm~\ref{alg:quNac}.   %To ensure the positive definitness of the preconditioner, 

The inputs are an initial point $x_0$, initial estimate $H_0$ and {\tt max\_q}; the maximum number of columns  allowed in  $\mathcal{S}_k$ at
each iteration
$k$.  In the first iteration, $k=0$, the search direction $d_0 =-H_0\nabla f_0$ is used. To determine $x_{k+1}$, a line search is used that first checks to see
if $a_k=1$ meets the line search criteria. In our implementation we use a sufficient descent criteria
\begin{equation} \label{eq:lnsch} 
f(x_k +a_k d_k) -f(x_k) \leq c_1 \alpha_k d_k^T \nabla f_k,
\end{equation}
with $c = 10^{-4}.$

 The PCG method~Algorithm~\ref{alg:conjgrad} is then called
with $H_k$ as a preconditioner to approximately solve $\nabla^2 f_{k+1} d_{k+1} = -\nabla f_{k+1}$ with the number of
iterations capped by {\tt max\_q}.
 Further limiting the number of PCG iterations is a tolerance
  \[\mbox{{\tt PCG\_tol}}=\min\left\{0.01, \norm{\nabla f(x_{k+1})}^{1/2}\right\},\] which
corresponds to the ``super-linear'' choice in inexact Newton methods~\cite{Dembo1982}.
The conjugate directions calculated during the PCG execution, which we denote by
$[p_{q(k)}, \ldots , p_{q(k+1)-1}]$ henceforth, are saved to form the columns of $\mathcal{S}_k.$
Specifically, the
 columns of $\mathcal{S}_k$ are the $\He_{k+1}-$normalized conjugate directions
 \begin{equation}\label{eq:Dscale}
 \mathcal{S}_k = \left[\frac{p_{q(k)}}{\norm{p_{q(k)}}_{\He_{k+1}}}, \ldots,
\frac{p_{q(k+1)-1}}{\norm{p_{q(k+1)-1}}_{\He_{k+1}}}\right].
\end{equation}
This normalization is done to simplify calculations, as with this
choice $\proj{\mathcal{S}_k}{\He_k} = \mathcal{S}_k \mathcal{S}_k^T.$  So that the resulting estimate is positive
definite, we only collect conjugate directions so long as
negative curvature is not encountered, line~\ref{ln:curv} of Algorithm~\ref{alg:conjgrad}. This ensures that $\mathcal{S}_k^T \He_{k+1}\mathcal{S}_k\succ 0.$ There is a safeguard for non-convex functions on line~\ref{ln:safeg} of Algorithm~\ref{alg:conjgrad}. If negative curvature is encountered on the first PCG iteration, then the first conjugate direction $p_0  = -H_k \nabla f_{k+1}$ is returned as the search direction.
Before moving onto the next iteration, the estimate matrix is updated  by either a full
or limited memory inverse quNac~\eqref{eq:invquNac} update, detailed in Sections~\ref{sec:qunacfull}
and~\ref{sec:qunaclim}, respectively.
%in Algorithms~\ref{alg:scaledquNacupdate} and~\ref{alg:LquNac}, respectively.

In line~\ref{ln:dAd} of Algorithm~\ref{alg:conjgrad}, we need to calculate a Hessian-vector product. This can be done efficiently through reverse AD
(\emph{Automatic Differentiation})~\cite{Christianson:1992}. Naturally there also exist problems and applications where
fast Hessian-vector products are readily available, such as Fast-Fourier transform, Neural
Networks~\cite{Pearlmutter1993a} or
obvious structure prevailing in the Hessian matrix. As a final option, the user would
be required to write an efficient subroutine for calculating Hessian-vector products.

\begin{algorithm}
\DontPrintSemicolon
\small
\caption{Newton-PCG quNac}
\label{alg:quNac}
\KwIn{  $H_0, x_0 \in \R,$ {\tt max\_q} $ \in \N$. }
$k =0,d_0 = -H_0 \nabla f_0$  \;
%$H_0 = C \cdot I$ \;
\While{$|\nabla f_k|/|\nabla f_0| > \epsilon$ or $|\nabla f_k|> \epsilon$ }{
Determine $a_k$ through a line-search on $\{a \, | \, x_k+ a d_k \}$ starting with $a_k =1$ \;
$x_{k+1} = x_k +a_k d_k$  \label{eq:stepsk}\;
$[\mathcal{S},\He_{k+1}\mathcal{S}, d_k] = $\texttt{PCG} $(\He_{k+1},H_k,x_{k+1},
\mbox{max\_$q$}, \mbox{{\tt PCG\_tol}})$ \label{eq:CGcall}\;
$H_{k+1} = \ref{eq:genquNac}(H_k,\He_{k+1}\mathcal{S} \rightarrow \mathcal{S} ) $, using~Algorithm~\ref{alg:scaledquNacupdate}\; \label{ln:Hkupdate}
$k = k+1$\;
}
\KwOut{$x_k$.}
\end{algorithm}

\begin{algorithm}
%\SetAlgoNoLine
\DontPrintSemicolon
\small
\caption{{\tt PCG}$(A,M^{-1},y_0, \mbox{max\_$q$}
,\mbox{{\tt PCG\_tol}})$}
\label{alg:conjgrad}
$r_0=\nabla f(y_0)$\;
$z_0= M^{-1}(r_0) $\;
$p_0 = -z_0$\;
$y_0 = 0$\;
\For{$i =0,\ldots, \mbox{max\_$q$}-1$ }{
${c}_i = \dotprod{A p_i ,p_i }$  \label{ln:dAd}\;
\If{\normalfont ${c}_i \leq 0$}{ \label{ln:curv}
   \lIf{$i>0$}{{\bf break} }
   \lElse{$y_0 = p_0$ \label{ln:safeg}}
}
%\lIf{\normalfont ${c}_i \leq 0$}{ }
$\displaystyle  \alpha_i = \frac{\dotprod{r_i ,z_i}}{{c}_i}$ \;
$y_{i+1} = y_i + \alpha_i p_i$\;
$\displaystyle r_{i+1} = r_i + \alpha_i Ap_i$\;
$z_{i+1} = M^{-1}r_{i+1}$ \;
$\displaystyle  \beta_i = \frac{\dotprod{r_{i+1} ,z_{i+1} }}{\dotprod{r_i ,z_i}}$ \label{alg:CGbeta}\;
$\displaystyle p_{i+1} =  -z_{i+1}  +  \beta_i p_i$\;
\If{$\norm{r_{i+1}}/\norm{r_0} < ${\tt PCG\_tol}}{
   $q = i+1$\;
	 {\bf break}\;
}
}
$q = \min\{ \mbox{max\_$q$},i \}$\;
\KwOut{$\mathcal{S} = \left[ {c}_0^{-1/2} p_0, \ldots, {c}_0^{-1/2} p_{q-1} \right], A\mathcal{S},y_{q} $.}
\end{algorithm}

\subsection{Full memory Inverse quNac} \label{sec:qunacfull}

Both the limited and full memory variants of the inverse quNac update have been implemented in a way that promotes
parallel linear algebra through Matrix multiplication.
To derive these two variants, let $\underline{\mathcal{S}}_k =\He_{k+1} \mathcal{S}_k$ be the $n\times q$ matrix
stored from executing PCG method in Algorithm~\ref{alg:conjgrad}.
With the normalization~\eqref{eq:Dscale} of $\mathcal{S}_k$, the inverse quNac update can be
calculated by
\begin{align*}
E_k & = \proj{\mathcal{S}_k}{ \He_{k+1}} + \proj{\mathcal{S}_k}{ \He_{k+1}}\He_{k+1}H_k( \He_{k+1} \proj{\mathcal{S}_k}{
\He_{k+1}} -I) -H_k \He_{k+1}
\proj{\mathcal{S}_k}{
\He_{k+1}} \\
&= \mathcal{S}_k\mathcal{S}_k^T +  \mathcal{S}_k \underline{\mathcal{S}}_k^{T}H_{k}( \underline{\mathcal{S}}_k
\mathcal{S}^T -I) -H_{k}
\underline{\mathcal{S}}_k \mathcal{S}_k^T \\
&= \mathcal{S}_k\left( I_{p \times p} +   \underline{\mathcal{S}}_k^{T} H_{k} \underline{\mathcal{S}}_k
\right)\mathcal{S}_k^T -H_{k}
\underline{\mathcal{S}}_k \mathcal{S}_k^T -
\mathcal{S}_k\underline{\mathcal{S}}_k^TH_{k}.
\end{align*}
This has been coded in Algorithm~\ref{alg:scaledquNacupdate} and costs $O(n^2 q)$ operations.
Line~\ref{ln:matrixHDproduct} is the bottleneck as it involves a multiplication of a possibly dense $n\times n$ matrix
with a $n \times q$ matrix. The cost of sequentially applying $q$ BFGS updates is also $O(n^2 q),$ the important
difference is that Algorithm~\ref{alg:scaledquNacupdate} can greatly benefit from multithreading and parallel
linear algebra, while there is no obvious parallelism in applying BFGS updates. In fact, if $q$ processors are available
in a shared memory architecture, then
the wall clock time of Algorithm~\ref{alg:scaledquNacupdate} is $O(n^2)$ plus additional overheads of the parallel
paradigm (such as creating and destroying threads).

\begin{algorithm}
 \DontPrintSemicolon
\small
\caption{Scaled Inverse $\ref{eq:genquNac}(H, \underline{\mathcal{S}} \rightarrow \mathcal{S} )$ update}
\label{alg:scaledquNacupdate}
\KwIn{  $H \in \R^{n \times n}$ and $\mathcal{S},\underline{\mathcal{S}} \in \R^{n \times q}$ }
  % $\underline{\mathcal{S}} = \He_k \mathcal{S}_k $\; %\quad \left(\mbox{ from PCG call}\right) $ \;
$\underline{H} = H\underline{\mathcal{S}}$ \label{ln:matrixHDproduct}\; % \quad \left(\mbox{costs $n^2 q$} \right)\\
$ \overline{H} =  \underline{H} \mathcal{S}^T$ \;  %\quad \left(\mbox{costs $n q^2$} \right)\\
$   E  = \mathcal{S}\left( I_{p \times p} +   \underline{\mathcal{S}}^{T} \underline{H} \right)\mathcal{S}^T
-\overline{H} -\overline{H}^T$\;
 % \quad \left(\mbox{costs $n q^2$} \right)
\KwOut{$E$.}
\end{algorithm}

The next Corollary shows that
when Algorithm~\ref{alg:quNac} uses \quNac updates, the resulting preconditioners satisfy the quadratic Hereditary property. Thus when Algorithm~\ref{alg:quNac}  is applied
to convex quadratic problems,  the method terminates after a total of $n$ inner steps of the PCG method. 

% The proof is essentially an application of the Restarting Preconditioner Lemma~\ref{lemma:restart}, %the quadratic Hereditary Corollary~\ref{cor:quadhered} and the comment after %Corollary~\ref{cor:quadhered}.

Due to this following Corollary, we chose to update the preconditioner with all available conjugate directions. This is in contrast with the strategies mentioned in~\cite{Morales2000a}, where 
the last conjugate directions or a uniform sampling of conjugate directions are used to perform L-BFGS updates.
\begin{corollary}[Quadratic Hereditary for quNac Preconditioner]\label{cor:metricgradconj}
Assume Algorithm~\ref{alg:quNac} is applied to a convex quadratic function $\phi(x)$ with $\nabla^2 \phi(x) \equiv Q \in \R^{n \times n}$,
and consider its $k$th major iteration, $k\geq 1$. Then $H_{k+1}Q\mathcal{S}_i =\mathcal{S}_i$ for $ i =0, \ldots, k.$
\end{corollary}
\proof We prove this using the Restarting Preconditioner Lemma~\ref{lemma:restart} to show that $\{p_0, \ldots,
p_{q(k+1)-1}\}$ is a $Q-$orthogonal set, then apply Corollary~\ref{cor:quadhered} and the comment after Corollary~\ref{cor:quadhered} to prove quadratic hereditary.
 The proof is by induction where our hypothesis is that the set $\{p_0, \ldots, p_{q(k)-1}\}$ is a $Q-$orthogonal set
and $p_j^T \nabla \phi(x_{k}) =0$ for all $0 \leq  j \leq q(k-1)-1.$

The base case of our induction is $k=2.$ The set of vectors $\{p_0, \ldots, p_{q(1)-1}\}$ calculated by the first PCG call are $Q-$orthogonal by construction.
At iteration $k=1$, as $x_1 + d_1$ is the minimum of the quadratic $\phi(x)$ over $x_1 \oplus \K_{q(1)-1},$ the step parameter $a_1=1$ is accepted.
Therefore $x_2 =x_1 + d_1$, $\nabla \phi(x_{2}) \in \K_{q(1)-1}^{\perp}$ and $ p_j^T  \nabla \phi(x_{2}) =0,$ for $j=0, \ldots q(1)-1.$
This proves, together with the action constraint $  H_{1} Qp_j  =
p_j$  for $j =0, \ldots, q(1)-1,$
that $x_2$ and $H_1$ are a restarting point and a restarting preconditioner, respectively, and by
Lemma~\ref{lemma:restart} the set $\{p_0, \ldots,p_{q(1)-1}, p_{q(1)}, \ldots,
p_{q(2)-1}\}$ is $Q-$orthogonal. This concludes the proof of our induction hypothesis for $k=2.$

Suppose that $p_j^T \nabla \phi(x_{k}) =0$ for all $0 \leq  j \leq
q(k-1)-1$ and $\{p_0, \ldots, p_{q(k)-1}\}$ are $Q-$orthogonal.
This $Q-$orthogonality guarantees by Corollary~\ref{cor:quadhered} that $H_{k}$  satisfies the hereditary property $
 Q H_{k} p_i = p_i$ for $i =0, \ldots, q(k)-1.$

At the $k$th iteration $a_k =1$ is accepted
as $x_{k} +d_k$  is the minimum of $x_k \oplus \, \mbox{span}\{p_{q(k-1)}, \ldots, p_{q(k)-1} \}$, thus
$p_{j}^T \nabla \phi(x_{k+1}) =0$ for $ q(k-1) \leq j<q(k).$ For $j<q(k-1)$ we have
\begin{align*}
p_j^T \nabla \phi(x_{k+1}) &= p_j^T\left( \nabla\phi(x_{k}) +  Q d_k \right) \\
 &= p_j^T\left( \nabla \phi(x_{k}) +  Q \left(\sum_{i=q(k-1)}^{q(k)-1}\alpha_i p_i\right)\right)\\
&=  p_j^T\nabla \phi(x_{k})+\sum_{i=q(k-1)}^{q(k)-1}\alpha_i\alpha_j p_j^T Q p_i  \quad \left(\mbox{applying the
induction
hypothesis}\right) \\
& = 0 + 0.
\end{align*}
Thus $x_{k+1}$ and $H_k$ are a restarting point and a restarting preconditioner, respectively, and
by Lemma~\ref{lemma:restart} the vectors
$\{p_{0} \ldots p_{q(k+1)-1}\}$ are $Q-$orthogonal, which concludes the induction. Finally, the columns of the sampling matrices are scalar multiples of the conjugate directions,
thus Corollary~\ref{cor:quadhered} and the comment that follow it guarantees the quadratic hereditary of $H_{k+1}$ is Algorithm~\ref{alg:quNac}.$\qed$

\subsection{Limited memory \quNac}\label{sec:qunaclim}

%Our limited memory method is essentially a parallel version of Morales and Nocedal's automatic L-BFGS
%preconditioner~\cite{Morales2000a}.

%It is also possible to adapt Algorithm~\ref{alg:quNac} for a limited memory implementation based on a L-BFGS
%preconditioner method~\cite{Morales2000a}.

To implement a limited memory variant of the
\texttt{inverse} quNac update~\eqref{eq:invquNac},
instead of updating $H_{k}$, in line~\ref{ln:Hkupdate} of Algorithm~\ref{alg:quNac},  we initiate $H_{k} =
H_0^{k+1}$ which is a user specified initial
estimate approximation (or simply the identity in the lack there of). Both $H_0^{k+1}$ and $H_{k+1}$ must be coded as
operators acting on vectors in $\R^n$ instead of explicit matrices.
In Algorithm~\ref{alg:LquNac} we show how to execute the operation $v \rightarrow H_0^{k+1}(v)+E_k(v)$ without
the need to store a matrix.
Let $\mathcal{S}_k$ and $\underline{\mathcal{S}}_k =\He_{k+1} \mathcal{S}_k$ be the $n\times q$ matrices
stored from the previous PCG call. Then to calculate $H_{k+1}v$ we have
\begin{align}
 (H_0^{k+1} +E_k)v & =  \proj{\mathcal{S}_k}{ \He_{k+1}} + (I-\proj{\mathcal{S}_k}{ \He_{k+1}}\He_{k+1})H_k(I- \He_{k+1}
\proj{\mathcal{S}_k}{
\He_{k+1}})\label{eq:LquNac} \\
&= \mathcal{S}_k\mathcal{S}_k^T v+  (I-\mathcal{S}_k \underline{\mathcal{S}}_k^{T})H_0^{k+1}(I-
\underline{\mathcal{S}}_k \mathcal{S}_k^T)v
\nonumber \\
&=   (H_0^{k+1}(v- \underline{\mathcal{S}}_k (\mathcal{S}_k^Tv))+\mathcal{S}_k \left( (\mathcal{S}_k^T v)-
\underline{\mathcal{S}}_k^{T}\left(H_0^{k+1}(v- \underline{\mathcal{S}}_k
(\mathcal{S}_k^Tv))\right)\right), \nonumber
\end{align}
which can be calculated efficiently by Algorithm~\ref{alg:LquNac}.
As the columns of $\mathcal{S}_k$ are $\He_{k+1}-$orthogonal, Proposition~\ref{prop:unravel} proves that
Algorithm~\ref{alg:LquNac}  has the same result, in exact precision, as applying the L-BFGS two-loop
recursion~\cite{Nocedal1980} to the columns of $\mathcal{S}_k$ and $\underline{\mathcal{S}}_k$. To compare the two
methods for
applying a preconditioner operator, we have placed the L-BFGS two-loop recursion and LquNac side-by-side in
Figure~\ref{fig:L}. The only difference between them is that  $\boxed{v}$ and $\boxed{r}$ in
Algorithm~\ref{alg:L-BFGS} are replaced by a new variable $z$ in Algorithm~\ref{alg:LquNac}. This small change
removes the dependency between the two lines in each {\tt for} loop in Algorithm~\ref{alg:L-BFGS} so that the loops can
be
calculated as matrix-vector products instead. Matrix-vector multiplications can be easily sped up
through multithreading and shared memory parallelism, while the two {\tt for} loops in Algorithm~\ref{alg:L-BFGS}
 are essentially sequential.

\begin{figure} \caption{Comparing the L-BFGS two-loop recursion with the parallel LquNac.} \label{fig:L}
\begin{flushleft}
 {\bf Input} $H_0:\R^n \rightarrow \R^n, \mathcal{S} =\left[s_1, \ldots, s_q \right], \underline{\mathcal{S}} = \left[
\underline{s}_1\ldots\underline{s}_q \right]$ and
$v \in \R^n.$\end{flushleft} \vspace{-0.3cm}
\begin{minipage}[t]{0.5\textwidth}
\begin{algorithm}[H]
\caption{ LquNac  }
\label{alg:LquNac}
 \DontPrintSemicolon \;  \PrintSemicolon
 $v^{\mathcal{S}} \leftarrow \mathcal{S}^Tv $\;
$z \leftarrow v - \underline{\mathcal{S}} v^{\mathcal{S}}$\;
%$\displaystyle z = \sum_{i=1}^{q} \underline{d}_i v^\mathcal{S}_i$\;
$\displaystyle  r \leftarrow H_0(z)$\;
 \DontPrintSemicolon \;   \PrintSemicolon
$r^{\underline{\mathcal{S}}} \leftarrow \underline{\mathcal{S}}^Tr $\;
$z \leftarrow r+  \mathcal{S}( v^{\mathcal{S}} -r^{\underline{\mathcal{S}}})$ \; \vspace{-0.1cm}
\KwOut{$z$}
\end{algorithm}
\end{minipage} %\hspace{1cm}
\begin{minipage}[t]{0.5\textwidth}
 \begin{algorithm}[H]
\caption{two-loop recursion}
\label{alg:L-BFGS}
%\KwIn{$H_0:\R^n \rightarrow \R^n, D =\left[s_1, \ldots, s_q \right], \underline{\mathcal{S}} = \left[ \underline{d}_1
\ldots
%\underline{d}_q \right]$ and $v \in \R^n$}
\For{$i = 1, \ldots, q$}{
$v^\mathcal{S}_i \leftarrow  s_i^T v$\;
$\boxed{v} \leftarrow v- v^\mathcal{S}_i \underline{d}_i$\;
}
$\displaystyle  r \leftarrow H_0(v)$\;
\For{$i = q, \ldots, 1$}{
$r^{\underline{\mathcal{S}}} \leftarrow  \underline{d}_i^T r$\;
$\boxed{r} \leftarrow r + s_i( v^\mathcal{S}_i - r^{\underline{\mathcal{S}}})$\;
}
\KwOut{$r$}
\end{algorithm}
\end{minipage}
\end{figure}

As of MATLAB version 7.4 (R2007a), MATLAB automatically multithreads matrix-vector multiplication, and tests on our
quad-core Desktop comparing the time taken to perform a L-BFGS two-loop recursion as compared to the LquNac update
 revealed that the speed-up can be more than four fold when there is sufficient number of columns in $\mathcal{S}_k$ and
$\underline{\mathcal{S}}_k$,
see
Figure~\ref{fig:ratioL}. This speed is specially important as applying this L-BFGS preconditioner is the bottle-neck in
the PCG iteration. There are a number of outliers in Figure~\ref{fig:ratioL} that are difficult to
investigate as multithreading is performed implicitly.  To have finer control and better exploit this parallelism an
explicit parallel
paradigm needs to be implemented, something we leave for future work.
\begin{figure} \centering
 \includegraphics[width=7cm]{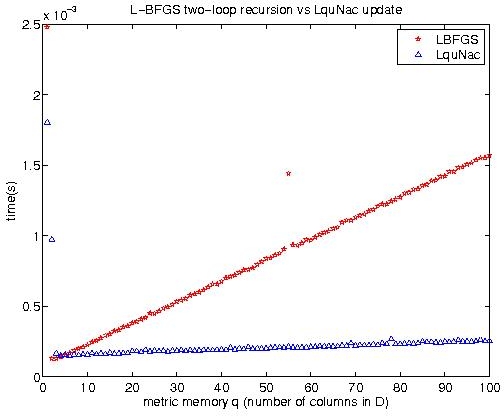}
\caption{Time taken by Applying the L-BFGS two-loop recursion in Algorithm~\ref{alg:L-BFGS} with the LquNac
update~\ref{alg:LquNac} where $\mathcal{S},\underline{\mathcal{S}}\in\R^{500 \times q}$ is randomly generated and $q$ is
increased from $1$
to $100$ and $H_0 =I.$ }
\label{fig:ratioL}
\end{figure}
Though we only consider this limited memory implementation that uses conjugate directions from the previous iteration,
certainly other  implementations are possible, for instance, by retaining conjugate directions from
other iterations.

\section{Numerical Tests} \label{sec:tests}

In our tests we compare five methods. The first two
methods are the full and limited memory \texttt{inverse} quNac update detailed in Algorithms~\ref{alg:quNac}. We
have labelled the two quNac methods by InverseQuNac and
InverseLQuNac, when the full memory variant in  Algorithm~\ref{alg:scaledquNacupdate}  and the limited variant in
Algorithm~\ref{alg:LquNac} are used to update the estimate, respectively. 
The third method is Newton\_CG implemented according to Algorithm 6.1 of~\cite{Nocedal1999}
though with an additional maximum number of CG iterations set to the dimension $n$ of the problem. The last two approaches are the BFGS and L-BFGS~\cite{Nocedal1980} methods. To compare the methods, we embed them in the same line search framework with a sufficient
descent criteria~\eqref{eq:lnsch} that initially checks if $a_k =1$ can be accepted. Though a line search that guarantees the Wolfe conditions is often advised for quasi-Newton methods, we found this to be inefficient when
applied to non-convex functions, as an almost exhaustive search for correct parameter $a_k$ would often occur. 
 The initial Hessian approximation was set to
\[ H_0 = \frac{\nabla f_0 ^T\nabla f_0}{\nabla f_0^T \He_0 \nabla f_0 }I.\]
In all the limited memory methods the maximum memory, {\tt max\_q} in the \quNac methods, was set to 20.

 Our MATLAB implementation ``quNac'' can be downloaded from the Edinburgh Research Group in
Optimization website: \url{http://www.maths.ed.ac.uk/ERGO}.
 In this package one can test different line search criteria, including Wolfe-conditions, and different initial Hessian $H_0$ approximations.

We have run tests on a Desktop with 64bit quad-core Intel(R) Core(TM) i5-2400S CPU @ 2.50GHz with 6MB cache size with a Scientific Linux
release 6.4 (Carbon) operating system.

%  A number of test results are reported below with our aforementioned 
% 
%  in particular the ``Classic Academic functions" were sensitive to
% the choice of line search, $H_0$ and the step method. We have selected the simplest results to report, but our settings
% should be considered a sample of the possible parameter space.

\subsection{Linear SVM with logistic loss}

Our first set of tests consists of convex Support Vector Machine (SVM) problems.
% Convex problems are a natural first
%choice, as they are common problems or sub-problems throughout optimization and have a unique %solution.
 SVMs have become a widely successful machine learning method for classification, and thanks to Chih-Chung Chang and Chih-Jen Lin LIBSVM
collection~\cite{Chang2011}, have readily
available data sets. We have selected all data sets for binary classification with less than or equal to $50'000$
features (dimensions).

The linear binary SVM problem consists of finding a separating hyperplane $f_w(x) = \dotprod{w,x}$ with $w \in \R^n$
that is able to predict the classification of $x\in X \subset \R^n$, namely, $f_w(x) >0$ and $f_w(x)\leq 0$ for the
first and second class, respectively. To this end, known data pairs $(x^i,y_i)$ are collected
where $x^i \in \R^n$ are \emph{feature} vectors and $y_i \in \{-1, 1\}$ are \emph{labels}, where $y_i$ indicates the
class of $x^i$ for $i =1, \ldots, m.$
The linear classifier $w$ is then selected based on these data pairs by minimizing a loss function, where a popular
choice~\cite{Yuan2012} is the logistic loss function
\begin{align*}
L_w(y,X) &= \sum_{i=1}^m \ln\left( 1 + \exp(-y_i \dotprod{x^i,w})\right).
\end{align*}
We use one of two regularizers, the $\ell_2$ norm
\[R_2(w) = \norm{w}_2^2,\]
or the \emph{pseudo-Huber} norm
\[R_{\mu}(w) = \mu \sum_{i=1}^n\left( \sqrt{1+\frac{x_i^2}{\mu^2}} -1 \right),\]
where $\mu<1.$ The pseudo-Huber norm is an approximation to the $\ell_1$ norm as $\mu \rightarrow 0,$ and has been shown
to be successful in promoting sparsity in convex regularized problems~\cite{Fountoulakis2013}. The resulting
unconstrained optimization problem is given by
\[ \min_w  L_w(y,X) + \lambda R_{\mu}(w), \]
where $\lambda$ is the regularizer parameter and has been set to $\lambda =1$ in all our tests.
 Our interest was in encountering the unique solution to these convex problems thus we solved the SVM problem with a
precision of $\epsilon =10^{-7}.$ We found through sampling a number of the problems that when increasing the
precision, the solution would become increasingly sparse up to approximately $\epsilon =10^{-7}$. Though
optimizing to a high tolerance raises the question of over-fitting, this is not an issue
here as the number of data points far exceeds the number of unknowns features, with the exception of the problem
\texttt{colon-cancer} ($62$ data points and $2000$ features) and {\tt duke breast cancer} ($44$ data points and $7129$ features).

In Tables~\ref{tab:LRe8} and~\ref{tab:LRe8Huber} we have the run times of each method to reach the unique solution
with a $\ell_2$ and pseudo-Huber regularizer, respectively.
In each table, ``ss'' represents ``small step'', in that the method takes
steps smaller than $\epsilon^2 = 10^{-14}$ before reaching the solution. While ``TO'' represents ``Timeout'' in that
the method exceeded the maximum time allowed, which we set to 10min. Each row corresponds to a problem and the
highlighted cells in the row indicate the smallest run time among all methods, while the boxed cell is the fastest  among
the limited memory methods. The last rows contain the standard deviation and average for each method across
all solved problems, though as each method failed to solve a number of problems, these statistics have to be interpreted with
care.

On the $\ell_2$ and pseudo-Huber regularized problems, InverseQuNac was the fastest method on most of the problems.
 Among the limited memory implementations, when tested on the $\ell_2$ regularized problems of Table~\ref{tab:LRe8}, Newton-CG was the fastest on 23,
InverseLQuNac was the fastest on 5 and L-BFGS was the fastest on 16 of the 44 problems tested.
 Though InverseLQuNac was the most
robust, failing to converge on only one problem and with the lowest standard deviation and average.
 For the pseudo-Huber regularized problems of Table~\ref{tab:LRe8Huber} the Newton-CG, InverseLQuNac and L-BFGS had the smallest run time on 11, 12
and 20 of the total 44 problems, respectively. The InverseLQuNac
was the robust out of the limited memory methods, failing only to converge on 3 problems, while Newton-CG and L-BFGS failed on 8 and 6 problems,
respectively.

 With the pseudo-Huber regularizer, as the sparse solution is approached, the Hessian becomes ill-conditioned~\cite{Fountoulakis2013}. This affected the
stability of Newton\_CG method. The InverseQuNac and InverseLQuNac seemed to be the least affected by this
ill-conditioning.

\begin{table} \centering	\tiny
 \begin{tabular}{|ccc|ccccc|}	\hline
		&	\# features	&	\# data	&		InverseQuNac	&
inverseLQuNac		&			Newton\_CG		&		BFGS	&
LBFGS		\\ \hline \hline
problem		&		&		&		Time(s)	&			Time(s)		&
	Time(s)		&		Time(s)	&			Time(s)		\\ \hline
a1a		&	119	&	1605	&		0.90	&			0.22		&
\cellcolor{blue!20}	\fbox{	0.17	}	&		1.74	&			0.38		\\
a2a		&	119	&	2265	&	\cellcolor{blue!20}	0.14	&			0.24
&		\fbox{	0.19	}	&		2.07	&			0.48		\\
a3a		&	122	&	3185	&	\cellcolor{blue!20}	0.16	&			0.31
&		\fbox{	0.27	}	&		2.69	&			0.58		\\
a4a		&	122	&	4781	&	\cellcolor{blue!20}	0.18	&			0.43
&		\fbox{	0.33	}	&		3.12	&			0.90		\\
a5a		&	122	&	6414	&	\cellcolor{blue!20}	0.25	&			0.52
&		\fbox{	0.45	}	&		4.20	&			1.09		\\
a6a		&	122	&	11220	&	\cellcolor{blue!20}	0.41	&			0.87
&		\fbox{	0.72	}	&		6.58	&			2.10		\\
a7a		&	122	&	16100	&	\cellcolor{blue!20}	0.60	&			1.32
&		\fbox{	1.23	}	&		9.71	&			3.49		\\
a8a		&	123	&	22696	&	\cellcolor{blue!20}	0.86	&			2.57
&		\fbox{	2.00	}	&		14.36	&			5.56		\\
a9a		&	123	&	32561	&	\cellcolor{blue!20}	1.31	&			4.13
&		\fbox{	3.46	}	&		21.89	&			9.48		\\
australian		&	14	&	690	&	\cellcolor{blue!20}	0.08	&
0.14		&		\fbox{	0.10	}	&		0.75	&			1.00
\\
australiansc		&	14	&	690	&	\cellcolor{blue!20}	0.05	&
0.07		&		\fbox{	0.06	}	&		0.21	&			0.12
\\
breast-cancer		&	10	&	683	&		0.02	&			0.02
&		\fbox{	0.02	}	&	\cellcolor{blue!20}	0.01	&			0.05
\\
breast-cancersc		&	10	&	683	&		0.12	&			0.17
&			0.15		&		0.20	&	\cellcolor{blue!20}	\fbox{	0.08	}
\\
cod-rna		&	8	&	59535	&	\cellcolor{blue!20}	0.91	&		\fbox{	1.63
}	&			1.99		&		8.09	&			8.73		\\
cod-rna.r		&	8	&	157413	&	\cellcolor{blue!20}	2.80	&		\fbox{
4.44	}	&			4.66		&		20.17	&			16.26
\\
colon-cancer		&	2000	&	62	&		1.65	&			0.24
&			0.26		&		42.68	&	\cellcolor{blue!20}	\fbox{	0.23	}
\\
covtype.binary		&	54	&	581012	&		10.38	&			16.36
&			20.56		&	\cellcolor{blue!20}	2.24	&		\fbox{	9.70	}
\\
covtype.binarysc		&	54	&	581012	&		12.22	&			19.83
&			19.56		&		35.45	&	\cellcolor{blue!20}	\fbox{	9.25	}
\\
diabetes		&	8	&	768	&	\cellcolor{blue!20}	0.03	&		\fbox{
0.04	}	&			0.32		&		0.20	&			0.18
\\
diabetessc		&	8	&	768	&	\cellcolor{blue!20}	0.03	&
0.04		&		\fbox{	0.04	}	&		0.13	&			0.05
\\
fourclass		&	2	&	862	&		0.02	&			0.03
&	\cellcolor{blue!20}	\fbox{	0.02	}	&		0.04	&			0.03
\\
fourclasssc		&	2	&	862	&	\cellcolor{blue!20}	0.02	&
0.02		&		\fbox{	0.02	}	&		0.03	&			0.02
\\
german.numer		&	24	&	1000	&	\cellcolor{blue!20}	0.06	&
0.12		&		\fbox{	0.12	}	&		0.99	&			2.31
\\
german.numersc		&	24	&	1000	&	\cellcolor{blue!20}	0.04	&
0.07		&		\fbox{	0.06	}	&		0.40	&			0.13
\\
gisettesc		&	5000	&	6000	&	\cellcolor{blue!20}	84.31	&		\fbox{
146.27	}	&			214.69		&		TO	&			161.39
\\
heart		&	13	&	270	&		0.07	&			0.08		&
\cellcolor{blue!20}	\fbox{	0.06	}	&		0.51	&			168.18		\\
heartsc		&	13	&	270	&	\cellcolor{blue!20}	0.02	&			0.04
&		\fbox{	0.04	}	&		0.15	&			0.06		\\
ionospheresc		&	34	&	351	&	\cellcolor{blue!20}	0.04	&
0.07		&		\fbox{	0.06	}	&		0.34	&			0.13
\\
liver-disorders		&	6	&	345	&	\cellcolor{blue!20}	0.05	&
0.07		&			0.06		&		0.08	&		\fbox{	0.05	}
\\
liver-disorderssc		&	6	&	345	&		0.04	&			0.07
&			0.06		&		0.11	&	\cellcolor{blue!20}	\fbox{	0.03	}
\\
mushrooms		&	112	&	8124	&		0.18	&			0.24
&			0.24		&		0.76	&	\cellcolor{blue!20}	\fbox{	0.17	}
\\
sonarsc		&	60	&	208	&	\cellcolor{blue!20}	0.04	&			0.08
&		\fbox{	0.07	}	&		0.25	&			0.61		\\
splice		&	60	&	1000	&	\cellcolor{blue!20}	0.05	&			0.09
&		\fbox{	0.09	}	&		0.46	&			ss		\\
splicesc		&	60	&	1000	&	\cellcolor{blue!20}	0.04	&		\fbox{
0.06	}	&			0.06		&		0.13	&			0.06
\\
svmguide1		&	4	&	3089	&		TO	&			TO
&			TO		&	\cellcolor{blue!20}	0.09	&		\fbox{	0.10	}
\\
svmguide3		&	22	&	1243	&	\cellcolor{blue!20}	0.04	&
0.07		&		\fbox{	0.06	}	&		0.40	&			0.21
\\
w1a		&	300	&	2477	&		0.16	&			0.20		&
	0.14		&		1.79	&	\cellcolor{blue!20}	\fbox{	0.13	}	\\
w2a		&	300	&	3470	&	\cellcolor{blue!20}	0.17	&			0.25
&			0.20		&		2.28	&		\fbox{	0.17	}	\\
w3a		&	300	&	4912	&	\cellcolor{blue!20}	0.21	&			0.28
&			0.27		&		2.47	&		\fbox{	0.24	}	\\
w4a		&	300	&	7366	&	\cellcolor{blue!20}	0.25	&			0.37
&			0.34		&		3.14	&		\fbox{	0.32	}	\\
w5a		&	300	&	9888	&	\cellcolor{blue!20}	0.29	&			0.48
&			0.46		&		3.76	&		\fbox{	0.41	}	\\
w6a		&	300	&	17188	&	\cellcolor{blue!20}	0.54	&			0.89
&			0.77		&		5.87	&		\fbox{	0.73	}	\\
w7a		&	300	&	24692	&	\cellcolor{blue!20}	0.78	&			1.28
&			1.44		&		8.75	&		\fbox{	1.12	}	\\
w8a		&	300	&	49749	&	\cellcolor{blue!20}	1.73	&			3.10
&			3.50		&		19.74	&		\fbox{	2.73	}	\\ \hline
 \multicolumn{3}{c}{ standard deviation}						&		12.94	&
\fbox{	22.42	}	&			32.78		&	\cellcolor{blue!20}	9.44	&
34.89		 \\
 \multicolumn{3}{c}{ average}						&	\cellcolor{blue!20}	2.84	&
\fbox{	4.83	}	&			6.50		&		5.33	&			9.51
 \\ \hline
\end{tabular}
	 \caption{Binary classification with $\ell_2$ regularizer and $\epsilon =10^{-7}$ and memory= 20. TO = TimeOut and ss = small step. The highlighted cells contain the fastest run time, while the boxed cells contain the
fastest run time among the limited memory implementations} \label{tab:LRe8}
\end{table}

\begin{table} \centering	\tiny
\begin{tabular}{|ccc|ccccc|}	\hline
		&	\# features	&	\# data	&		InverseQuNac	&
inverseLQuNac		&			Newton\_CG		&		BFGS	&
LBFGS		\\ \hline \hline
problem		&		&		&		Time(s)	&			Time(s)		&
	Time(s)		&		Time(s)	&			Time(s)		\\ \hline
a1a		&	119	&	1605	&	\cellcolor{blue!20}	3.10	&		\fbox{	15.98
}	&			33.44		&		6.38	&			ss		\\
a2a		&	119	&	2265	&	\cellcolor{blue!20}	2.89	&		\fbox{	12.34
}	&			54.95		&		6.41	&			ss		\\
a3a		&	122	&	3185	&	\cellcolor{blue!20}	4.22	&			13.60
&			119.53		&		6.47	&		\fbox{	7.02	}	\\
a4a		&	122	&	4781	&	\cellcolor{blue!20}	4.03	&			38.26
&			176.66		&		9.05	&		\fbox{	6.46	}	\\
a5a		&	122	&	6414	&	\cellcolor{blue!20}	3.96	&			16.34
&			77.39		&		11.20	&		\fbox{	7.43	}	\\
a6a		&	122	&	11220	&	\cellcolor{blue!20}	6.45	&			18.89
&			118.58		&		17.13	&		\fbox{	10.00	}	\\
a7a		&	122	&	16100	&	\cellcolor{blue!20}	7.82	&			25.54
&			188.65		&		18.51	&		\fbox{	18.45	}	\\
a8a		&	123	&	22696	&	\cellcolor{blue!20}	8.26	&		\fbox{	20.89
}	&			TO		&		25.02	&			30.30		\\
a9a		&	123	&	32561	&	\cellcolor{blue!20}	12.14	&			27.44
&			TO		&		34.28	&		\fbox{	16.68	}	\\
australian		&	14	&	690	&	\cellcolor{blue!20}	0.10	&
0.14		&		\fbox{	0.12	}	&		0.80	&			0.93
\\
australiansc		&	14	&	690	&	\cellcolor{blue!20}	0.04	&
0.07		&		\fbox{	0.07	}	&		0.42	&			0.15
\\
breast-cancer		&	10	&	683	&		0.02	&			0.02
&		\fbox{	0.02	}	&	\cellcolor{blue!20}	0.01	&			0.05
\\
breast-cancersc		&	10	&	683	&		0.36	&			0.94
&			1.61		&		0.40	&	\cellcolor{blue!20}	\fbox{	0.23	}
\\
cod-rna		&	8	&	59535	&	\cellcolor{blue!20}	0.99	&		\fbox{	1.87
}	&			3.41		&		7.51	&			7.21		\\
cod-rna.r		&	8	&	157413	&	\cellcolor{blue!20}	2.35	&		\fbox{
3.96	}	&			5.05		&		17.08	&			13.44
\\
colon-cancer		&	2000	&	62	&		58.73	&	\cellcolor{blue!20}	\fbox{
26.77	}	&			319.28		&		261.38	&			436.45
\\
covtype.binary		&	54	&	581012	&		9.24	&			14.40
&			18.32		&	\cellcolor{blue!20}	1.95	&		\fbox{	8.57	}
\\
covtype.binarysc		&	54	&	581012	&		563.51	&			TO
&			TO		&		TO	&	\cellcolor{blue!20}	\fbox{	210.94	}
\\
diabetes		&	8	&	768	&	\cellcolor{blue!20}	0.04	&		\fbox{
0.05	}	&			0.36		&		0.21	&			0.19
\\
diabetessc		&	8	&	768	&	\cellcolor{blue!20}	0.06	&		\fbox{
0.10	}	&			0.12		&		0.25	&			0.12
\\
fourclass		&	2	&	862	&		0.03	&			0.03
&	\cellcolor{blue!20}	\fbox{	0.02	}	&		0.03	&			0.03
\\
fourclasssc		&	2	&	862	&	\cellcolor{blue!20}	0.02	&
0.03		&			0.03		&		0.05	&		\fbox{	0.03	}
\\
german.numer		&	24	&	1000	&	\cellcolor{blue!20}	0.08	&
0.18		&		\fbox{	0.16	}	&		1.02	&			2.48
\\
german.numersc		&	24	&	1000	&	\cellcolor{blue!20}	0.10	&
0.16		&		\fbox{	0.15	}	&		0.64	&			0.23
\\
gisettesc		&	5000	&	6000	&		TO	&			TO
&			TO		&		TO	&			TO		\\
heart		&	13	&	270	&	\cellcolor{blue!20}	0.07	&			0.09
&		\fbox{	0.08	}	&		0.57	&			0.82		\\
heartsc		&	13	&	270	&	\cellcolor{blue!20}	0.08	&			0.16
&		\fbox{	0.15	}	&		0.36	&			0.19		\\
ionospheresc		&	34	&	351	&	\cellcolor{blue!20}	0.26	&		\fbox{
0.67	}	&			2.29		&		1.10	&			ss
\\
liver-disorders		&	6	&	345	&		0.16	&			0.52
&			0.64		&		0.32	&	\cellcolor{blue!20}	\fbox{	0.15	}
\\
liver-disorderssc		&	6	&	345	&		0.20	&			1.50
&			0.88		&		0.33	&	\cellcolor{blue!20}	\fbox{	0.11	}
\\
mushrooms		&	112	&	8124	&		11.88	&			27.17
&			284.39		&		11.24	&	\cellcolor{blue!20}	\fbox{	5.36	}
\\
sonarsc		&	60	&	208	&	\cellcolor{blue!20}	0.80	&		\fbox{	4.54
}	&			5.94		&		ss	&			ss		\\
splice		&	60	&	1000	&	\cellcolor{blue!20}	0.13	&			0.24
&		\fbox{	0.18	}	&		0.62	&			0.39		\\
splicesc		&	60	&	1000	&	\cellcolor{blue!20}	0.11	&
0.20		&		\fbox{	0.18	}	&		0.47	&			0.18
\\
svmguide1		&	4	&	3089	&		TO	&			TO
&			TO		&		0.42	&	\cellcolor{blue!20}	\fbox{	0.17	}
\\
svmguide3		&	22	&	1243	&	\cellcolor{blue!20}	0.78	&
127.34		&		\fbox{	5.86	}	&		1.79	&			ss
\\
w1a		&	300	&	2477	&	\cellcolor{blue!20}	9.62	&		\fbox{	30.97
}	&			469.14		&		23.72	&			46.92		\\
w2a		&	300	&	3470	&	\cellcolor{blue!20}	10.30	&			26.77
&			236.00		&		26.79	&		\fbox{	24.87	}	\\
w3a		&	300	&	4912	&	\cellcolor{blue!20}	15.43	&			52.27
&			458.37		&		32.16	&		\fbox{	25.74	}	\\
w4a		&	300	&	7366	&	\cellcolor{blue!20}	18.99	&			58.65
&			189.74		&		42.56	&		\fbox{	47.68	}	\\
w5a		&	300	&	9888	&	\cellcolor{blue!20}	23.30	&			44.99
&			TO		&		32.60	&		\fbox{	32.83	}	\\
w6a		&	300	&	17188	&	\cellcolor{blue!20}	23.28	&		\fbox{	38.59
}	&			355.64		&		48.79	&			46.40		\\
w7a		&	300	&	24692	&	\cellcolor{blue!20}	28.32	&			82.81
&			TO		&		81.88	&		\fbox{	61.41	}	\\
w8a		&	300	&	49749	&	\cellcolor{blue!20}	61.45	&			124.26
&			TO		&		147.90	&		\fbox{	74.72	}	\\ \hline
 \multicolumn{3}{c}{ standard deviation}						&		86.82	&
\cellcolor{blue!20}	\fbox{	30.85	}	&			137.86		&		47.13	&
\multicolumn{1}{c}{			77.27		} \\
 \multicolumn{3}{c}{ average}						&		21.28	&
\cellcolor{blue!20}	\fbox{	20.97	}	&			86.87		&		21.46	&
\multicolumn{1}{c}{			30.14		} \\ \hline
\end{tabular}	
 \caption{Binary classification with pseudo-Huber regularizer and $\epsilon =10^{-7}$ and memory$= 20$. TO = TimeOut and ss = small step. The highlighted cells contain the fastest run time, while the boxed cells contain the
fastest run time among the limited memory implementations} \label{tab:LRe8Huber}
\end{table}

To appraise the rate of convergence of each method,  in Figure~\ref{fig:Huber_epsilon}  we have plotted the evolution of
the error through time for each method applied to the \texttt{epsilon\_normalized} problem.  The \texttt{epsilon\_normalized} problem is the most challenging of our SVM problems.
Originating from the Pascal Large Scale Learning Challenge 2008\footnote{\url{http://largescale.ml.tu-berlin.de/about/}}, \texttt{epsilon\_normalized} is very ill-conditioned. The
L-BFGS and InversequNac enjoyed the fastest convergence, though the L-BFGS method suffered from some oscillation thus the quality of its solution depends on when the algorithm is terminated.

In Figure~\ref{fig:L_2_cod-rna} we have plotted the evolution of
the error through time for the full memory methods: InverseQuNac, BFGS, and Newton-CG, applied to {\tt cod-rna.r} with an $\ell_2$ regularizer. In this plot, the InverseQuNac method converges first in
just over $2$ seconds followed by  Newton\_CG in $4$ seconds. The BFGS method needs more than $16$ seconds to converge.

% REDO TEST with Larger markers, needs to be clear in black and white -->
To not forget the benefits of limited memory implementations, we have tested two additional large-scale problems,
 \texttt{rcv1\_train-binary} and \texttt{duke breast-cancer}, whose dimensions do not permit
a full memory implementation.
In Figures~\ref{fig:L_2rcv1_train} and~\ref{fig:Huber_duke} we have plotted the evolution of
the error through time for InverseLQuNac, Newton\_CG and L-BFGS. %

 The three methods had similar
results on the \texttt{rcv1\_train-binary} though the L-BFGS converged first.
While on the \texttt{duke breast-cancer}, the InverseLQuNac converged in just over 60 seconds, Newton-CG
stagnated at a very high error of 0.4 and L-BFGS rapidly decreased the error initially, but stagnated at 
an error of $10^{-6}.$

\begin{figure} \centering
\begin{tikzpicture}
\node (fig) at (0,0) {\includegraphics[width =  \textwidth]{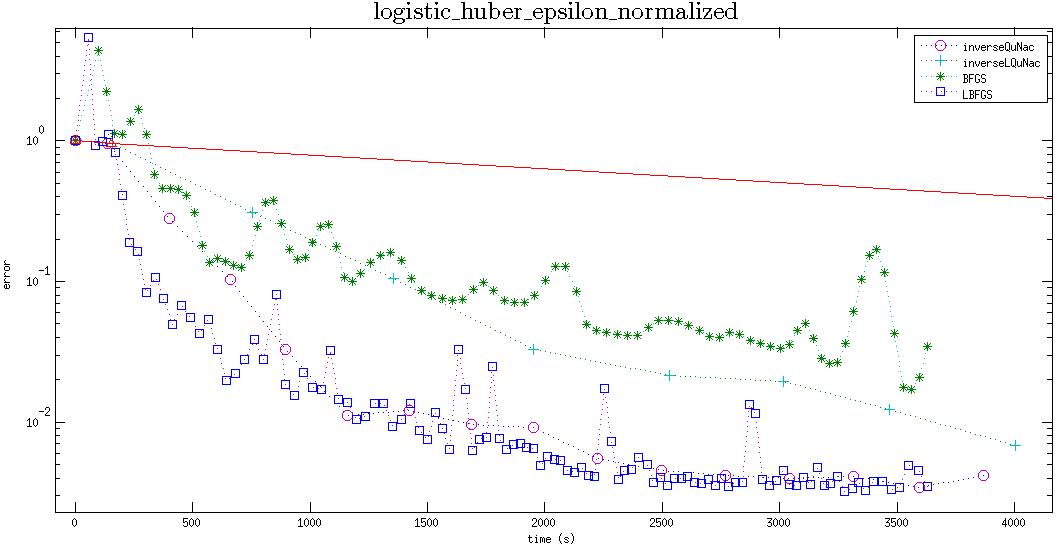}};
\node (lab) at (2,1.7) {\tiny Newton\_CG};
\end{tikzpicture}
\caption{The \texttt{epsilon\_normalized} problem with pseudo-Huber regularizer has
 $400,000$ data points and $2000$ features.}
\label{fig:Huber_epsilon}
\end{figure}

\begin{figure}
	\centering
        \begin{subfigure}[t]{0.7\textwidth}
                \includegraphics[width=\textwidth]{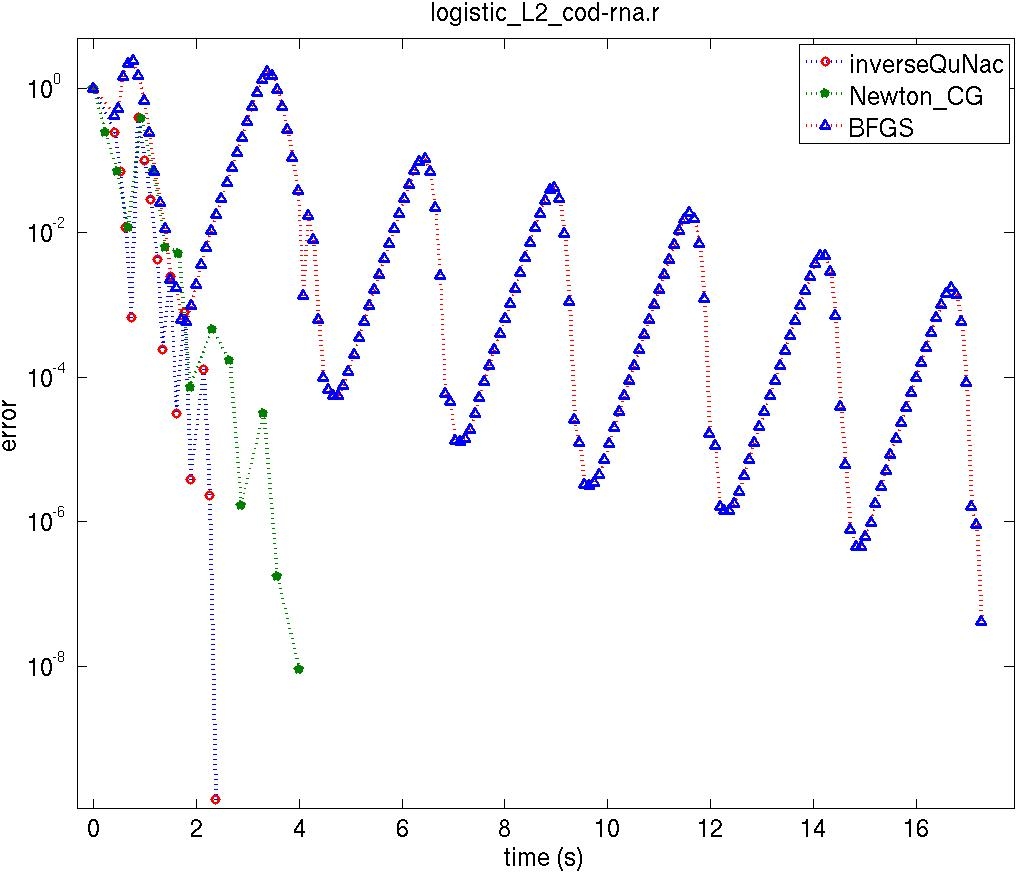}
\caption{The evolution of the error through time for each method applied to SVM with $\ell_2$
regularizer  on the \texttt{cod-rna.r2} problem. The error is on a logarithmic scale.}
                \label{fig:L_2_cod-rna}
        \end{subfigure} \\%
        \begin{subfigure}[t]{0.45\textwidth}
                \includegraphics[width=\textwidth]{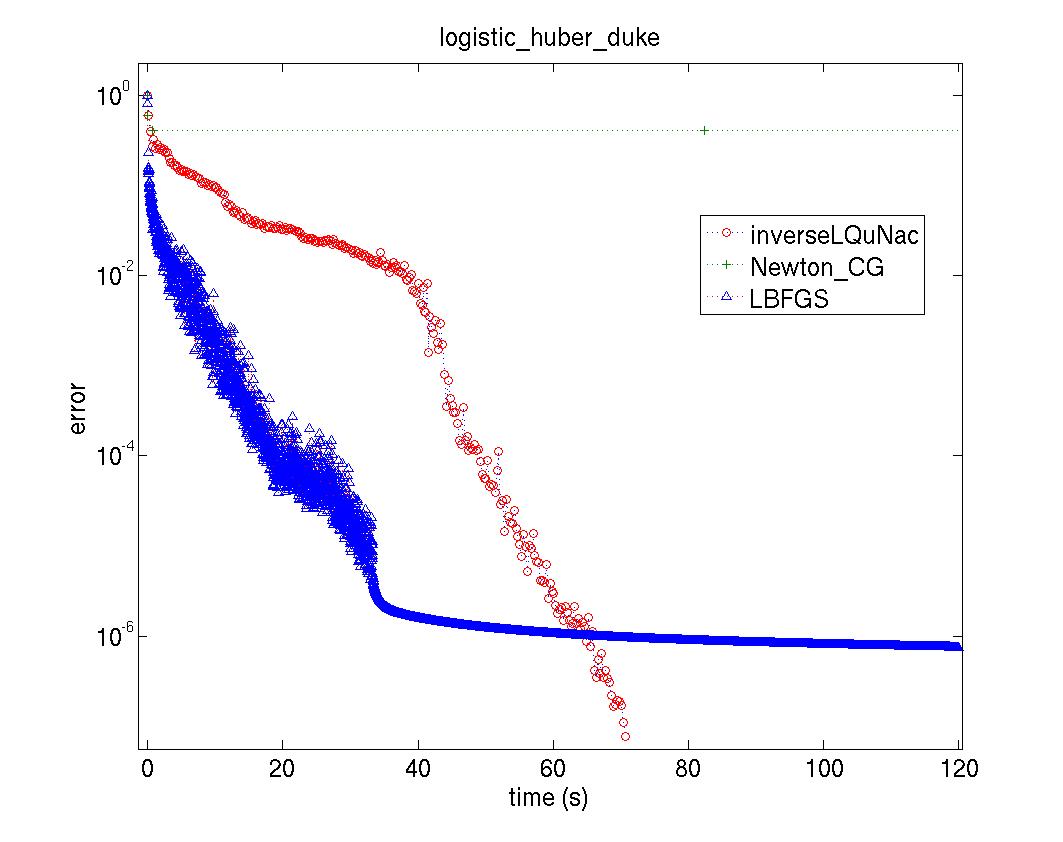}
                \caption{The \texttt{duke breast-cancer} problem with pseudo-Huber regularizer has $44$ data points and $7129$ features.}
                \label{fig:Huber_duke}
        \end{subfigure}%
				\hspace{0.5cm}
        \begin{subfigure}[t]{0.45\textwidth}
                \includegraphics[width=\textwidth]{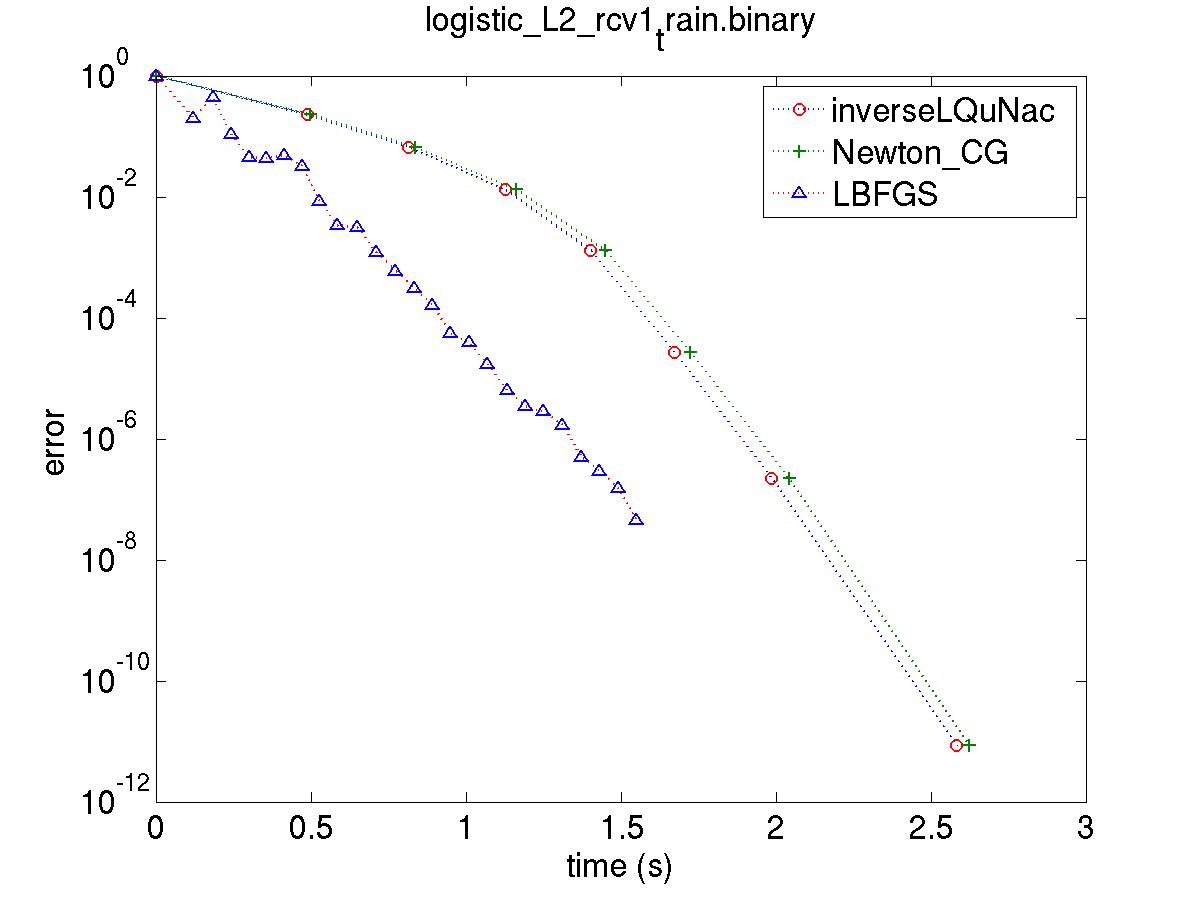}
                \caption{The \texttt{rcv1\_train-binary} problem with $\ell_2$ regularizer has $20242$ data points and
$47236$ features.}
                \label{fig:L_2rcv1_train}
        \end{subfigure}%
\caption{The evolution of error through time for each limited memory method applied to SVM LR problem}
\end{figure}

\subsection{Classic Academic functions}

We selected a number of academic unconstrained problems from~\cite{Mori1981} based solely on scalability of the function
and availability of the MATLAB code, in that, together  with their derivatives were
readily coded thanks to John Burkardt (\url{http://people.sc.fsu.edu/~jburkardt/m_src/test_opt/test_opt.html}),
see Table~\ref{tab:prodesc}.
Among these tests were two convex quadratic functions with ill-conditioned  $Q \in \R^{n\times n}$ Hessian matrices;
The Hilbert matrix $Q^H_{ij} = 2/ ( i + j - 1 )$ for $1 \leq i,j\leq n$ and the Gregory and Karney Tridiagonal Matrix
where $Q_{11} =4$, $Q_{12} =-2$, $Q_{ii} =2$, $Q_{i(i+1)} =-2 = Q_{i(i-1)}$ for $i=2, \ldots, n.$

\begin{table}
 \begin{tabular}{l|l} \hline
Problem	&	Description		\\ \hline
The Watson function	&	quartic function		\\
The Penalty Function \#1	&	quartic penalty function		\\
The Penalty Function \#2	&	nonlinear penalty function		\\
The Trigonometric Function	&	squared sum of trig. Functions		\\
The Extended Rosenbrock parabolic valley \#1	&	indefinite Hessian matrix		\\
The Extended Powell Singular Quartic	&	Singular Hessian matrix		\\
The Chebyquad Function	&	 quadrature of Chebyshev polynomials with no known solution		\\ \hline
The Gregory and Karney Tridiagonal Matrix	&	Ill-conditioned positive definite quadratic		\\
The Hilbert Matrix Function	&	Ill-conditioned positive definite quadratic
 \end{tabular}
\caption{Unconstrained test set description} \label{tab:prodesc}
\end{table}
Each test specifies an initial starting point from which we run each method until $\norm{\nabla f(x)}/
\norm{\nabla f(x_0)} < \epsilon$, which we set to $\epsilon = 10^{-8}$, or until 10 minutes of time was exhausted.
As a number of these problems were not convex, we employed a resetting and curvature criteria. Before taking a step in
the $d_k$ direction, line~\ref{eq:stepsk} of Algorithm~\ref{alg:quNac}, we verify if
\[ -\frac{\dotprod{d_k,\nabla f_k} }{\norm{d_k}\norm{\nabla f_k}} > \epsilon,\]
otherwise we reset the estimate $H_k = H_0$ and set $d_k = - H_0\nabla f_k$. As many of these test functions have
indefinite Hessian matrices, we terminate the PCG method at line~\ref{ln:curv} of Algorithm~\ref{alg:conjgrad} when
negative curvature
$\dotprod{Ap_i, p_i} <0$ is encountered. If no direction of positive curvature is encountered, the estimate matrix is
not updated, and we repeat the use of the previous estimate matrix $H_{k+1}=H_k.$ This idea of repeating a previous estimate
has been analysed in detail and tested in~\cite{Gill2001}.

In Table~\ref{tab:aca10-8} we report times taken to attain a stationary point for each method. The Newton\_CG method
was the fastest on 31 out of the 66 problems, while InverseQuNac, InverseLQuNac, BFGS and L-BFGS methods were the
fastest on 15, 5, 7 and 8 problems, respectively. Comparing only the limited memory methods, Newton-CG, InverseLQuNac
and L-BFGS methods were the fastest on 36, 20 and 8 problems, respectively. The InverseQuNac is the most stable, in that, it reached a stationary point on the largest number of problems; 65 out of 66. The results show that this particular adaptation of
the quNac method for general non-convex functions was very robust.

\begin{table} \centering	\tiny \vspace{-2.5cm}   \begin{tabular}{|cc|ccccc|}	 \hline
Problem		&	dimension	&		InverseQuNac	&			inverseLQuNac
&			Newton\_CG		&		BFGS	&			LBFGS		\\
\hline \hline
The Penalty Function \#2		&	100	&		0.16	&			0.18
&			0.24		&		0.59	&	\cellcolor{blue!20}	\fbox{	0.06	}
\\
		&	125	&		0.25	&			0.19		&
0.28		&		1.32	&	\cellcolor{blue!20}	\fbox{	0.09	}	\\
		&	150	&		0.33	&			0.26		&
0.38		&		2.35	&	\cellcolor{blue!20}	\fbox{	0.12	}	\\
The Penalty Function \#1		&	100	&		0.06	&			0.06
&	\cellcolor{blue!20}	\fbox{	0.05	}	&		0.06	&			0.05
\\
		&	200	&		0.08	&			0.06		&
\cellcolor{blue!20}	\fbox{	0.05	}	&		0.11	&			0.05		\\
		&	300	&		0.10	&			0.06		&
\cellcolor{blue!20}	\fbox{	0.05	}	&		0.17	&			0.05		\\
		&	400	&		0.12	&			0.06		&
\cellcolor{blue!20}	\fbox{	0.05	}	&		0.23	&			0.05		\\
		&	500	&		0.18	&			0.07		&
\cellcolor{blue!20}	\fbox{	0.05	}	&		0.32	&			0.05		\\
		&	600	&		0.22	&			0.06		&
0.05		&		0.39	&	\cellcolor{blue!20}	\fbox{	0.05	}	\\
		&	700	&		0.28	&			0.07		&
0.05		&		0.51	&	\cellcolor{blue!20}	\fbox{	0.05	}	\\
		&	800	&		0.35	&			0.07		&
0.05		&		0.64	&	\cellcolor{blue!20}	\fbox{	0.05	}	\\
		&	900	&		0.45	&			0.07		&
0.05		&		0.79	&	\cellcolor{blue!20}	\fbox{	0.05	}	\\
		&	1000	&		0.56	&			0.07		&
0.05		&		0.98	&	\cellcolor{blue!20}	\fbox{	0.05	}	\\
 Rosenbrock \# 1		&	100	&		0.07	&			0.09		&
\cellcolor{blue!20}	\fbox{	0.07	}	&		0.09	&			0.08		\\
		&	200	&		0.09	&			0.10		&
\cellcolor{blue!20}	\fbox{	0.07	}	&		0.19	&			0.09		\\
		&	300	&		0.13	&			0.11		&
\cellcolor{blue!20}	\fbox{	0.08	}	&		0.27	&			0.09		\\
		&	400	&		0.17	&			0.11		&
\cellcolor{blue!20}	\fbox{	0.08	}	&		0.39	&			0.10		\\
		&	500	&		0.25	&			0.12		&
\cellcolor{blue!20}	\fbox{	0.09	}	&		0.52	&			0.10		\\
		&	600	&		0.31	&			0.12		&
\cellcolor{blue!20}	\fbox{	0.09	}	&		0.68	&			0.11		\\
		&	700	&		0.40	&			0.12		&
\cellcolor{blue!20}	\fbox{	0.09	}	&		0.85	&			0.11		\\
		&	800	&		0.50	&			0.12		&
\cellcolor{blue!20}	\fbox{	0.10	}	&		1.05	&			0.11		\\
		&	900	&		0.63	&			0.13		&
\cellcolor{blue!20}	\fbox{	0.10	}	&		1.31	&			0.11		\\
		&	1000	&		0.77	&			0.13		&
\cellcolor{blue!20}	\fbox{	0.11	}	&		1.64	&			0.11		\\
The Extended Powell		&	100	&		0.08	&			0.08		&
\cellcolor{blue!20}	\fbox{	0.07	}	&		0.12	&			0.16		\\
		&	200	&		0.09	&			0.08		&
\cellcolor{blue!20}	\fbox{	0.08	}	&		0.28	&			0.11		\\
		&	300	&		0.12	&			0.09		&
\cellcolor{blue!20}	\fbox{	0.08	}	&		0.41	&			0.11		\\
		&	400	&		0.15	&			0.09		&
\cellcolor{blue!20}	\fbox{	0.09	}	&		0.58	&			0.30		\\
		&	500	&		0.20	&			0.10		&
\cellcolor{blue!20}	\fbox{	0.09	}	&		0.78	&			0.31		\\
		&	600	&		0.28	&	\cellcolor{blue!20}	\fbox{	0.10	}	&
	0.10		&		1.00	&			0.32		\\
		&	700	&		0.33	&	\cellcolor{blue!20}	\fbox{	0.10	}	&
	0.11		&		1.24	&			0.33		\\
		&	800	&		0.42	&	\cellcolor{blue!20}	\fbox{	0.10	}	&
	0.11		&		1.52	&			0.33		\\
		&	900	&		0.53	&	\cellcolor{blue!20}	\fbox{	0.11	}	&
	0.11		&		1.95	&			0.34		\\
		&	1000	&		0.67	&	\cellcolor{blue!20}	\fbox{	0.11	}	&
	0.12		&		2.45	&			0.36		\\
The Watson function		&	100	&		1.07	&		\fbox{	2.43	}	&
	6.27		&	\cellcolor{blue!20}	0.93	&			TO		\\
		&	200	&		7.73	&		\fbox{	15.35	}	&
20.31		&	\cellcolor{blue!20}	1.68	&			TO		\\
		&	300	&		9.43	&		\fbox{	60.54	}	&
63.27		&	\cellcolor{blue!20}	2.78	&			TO		\\
		&	400	&		65.36	&			95.45		&		\fbox{
74.20	}	&	\cellcolor{blue!20}	3.42	&			TO		\\
		&	500	&		97.53	&		\fbox{	311.24	}	&
344.94		&	\cellcolor{blue!20}	5.01	&			TO		\\
		&	600	&		71.71	&			328.34		&		\fbox{
163.18	}	&	\cellcolor{blue!20}	6.69	&			TO		\\
The Chebyquad Function		&	10	&		0.28	&		\fbox{	0.42	}	&
	0.45		&	\cellcolor{blue!20}	0.20	&			0.60		\\
		&	20	&	\cellcolor{blue!20}	0.15	&			0.77		&
\fbox{	0.74	}	&		ss	&			ss		\\
		&	30	&	\cellcolor{blue!20}	0.81	&			TO		&
\fbox{	23.21	}	&		ss	&			ss		\\
Tridiagonal Matrix Function		&	100	&		0.05	&			0.07
&	\cellcolor{blue!20}	\fbox{	0.02	}	&		ss	&			TO
\\
		&	200	&		0.08	&			0.14		&
\cellcolor{blue!20}	\fbox{	0.05	}	&		ss	&			TO		\\
		&	300	&		0.17	&			0.24		&
\cellcolor{blue!20}	\fbox{	0.07	}	&		ss	&			TO		\\
		&	400	&		0.27	&			0.35		&
\cellcolor{blue!20}	\fbox{	0.10	}	&		ss	&			TO		\\
		&	500	&		0.55	&			0.48		&
\cellcolor{blue!20}	\fbox{	0.13	}	&		ss	&			TO		\\
		&	600	&		0.75	&			0.57		&
\cellcolor{blue!20}	\fbox{	0.17	}	&		ss	&			TO		\\
		&	700	&		1.05	&			0.69		&
\cellcolor{blue!20}	\fbox{	0.20	}	&		ss	&			TO		\\
		&	800	&		1.42	&			0.82		&
\cellcolor{blue!20}	\fbox{	0.24	}	&		ss	&			TO		\\
		&	900	&		2.18	&			0.97		&
\cellcolor{blue!20}	\fbox{	0.27	}	&		ss	&			TO		\\
		&	1000	&		3.12	&			1.13		&
\cellcolor{blue!20}	\fbox{	0.31	}	&		ss	&			TO		\\
The Hilbert Matrix Function		&	100	&	\cellcolor{blue!20}	0.03	&		\fbox{
0.04	}	&			0.05		&		0.30	&			19.71
\\
		&	200	&	\cellcolor{blue!20}	0.07	&		\fbox{	0.08	}	&
	0.18		&		0.93	&			162.88		\\
		&	300	&	\cellcolor{blue!20}	0.12	&		\fbox{	0.29	}	&
	0.47		&		2.15	&			TO		\\
		&	400	&	\cellcolor{blue!20}	0.21	&		\fbox{	0.54	}	&
	0.69		&		3.46	&			TO		\\
		&	500	&	\cellcolor{blue!20}	0.34	&		\fbox{	0.83	}	&
	1.70		&		6.07	&			549.74		\\
		&	600	&	\cellcolor{blue!20}	0.46	&		\fbox{	1.31	}	&
	2.30		&		8.28	&			538.13		\\
		&	700	&	\cellcolor{blue!20}	0.60	&		\fbox{	1.78	}	&
	3.15		&		11.06	&			TO		\\
		&	800	&	\cellcolor{blue!20}	0.79	&		\fbox{	2.28	}	&
	4.53		&		14.16	&			TO		\\
		&	900	&	\cellcolor{blue!20}	1.03	&		\fbox{	2.97	}	&
	5.17		&		17.84	&			TO		\\
		&	1000	&	\cellcolor{blue!20}	1.22	&		\fbox{	3.53	}	&
	5.70		&		21.76	&			TO		\\
The Trigonoestimate Function		&	100	&		ss	&			2.44
&	\cellcolor{blue!20}	\fbox{	2.34	}	&		ss	&			ss
\\
		&	200	&	\cellcolor{blue!20}	0.61	&			ss		&
	ss		&		ss	&			ss		\\
		&	300	&	\cellcolor{blue!20}	1.50	&			ss		&
	ss		&		ss	&			ss		\\
		&	400	&	\cellcolor{blue!20}	1.70	&			23.45		&
\fbox{	17.36	}	&		ss	&			ss		\\ \hline
\multicolumn{2}{c}{standard deviation}				&		16.72	&			57.71
&		\fbox{	48.44	}	&	\cellcolor{blue!20}	4.55	& \multicolumn{1}{c}{
124.76		}\\ \hline
\end{tabular}  % \input{tab_aca10-8}
	 \caption{\small Tests on Academic functions from Table~\ref{tab:prodesc} with $\epsilon = 10^{-8}$ and memory$=20$.  TO = TimeOut and ss = small step. The highlighted cells contain the fastest run time, while the boxed cells contain the fastest run time among the limited memory implementations } \label{tab:aca10-8}
\end{table}

% Unconstrained methods such as the quNac methods are often embedded in globalizing strategies such as an
% Interior-Point,
% Trust-Region or Sequential-Quadratic programming framework.
% It is in these setting where the true use of a method quNac can be tested. Though such embedded tests present
% challenges: How much of the success of a such a combined implementation is due to the unconstrained solver?

\section{Conclusion} \label{sec:conclusion}

We have developed a family of updating schemes that generates a sequence of symmetric matrices which approximate a desired target sequence of symmetric matrices, where only the action of our target matrices on certain subspaces is known.
Furthermore, the updates have small rank, with rank at most three times that of the given subspace dimension.
This setup allows us to estimate the inverse of a matrix field, such as the inverse Hessian matrix, only by sampling its
action and never explicitly calculating the inverse. Sufficient conditions for positive definiteness and the quadratic hereditary property of the estimates are established in this general setting.

The application we focus on is solving sequences of Newton systems; a common building block of many optimization methods. In this setting, we match the action of our estimate matrix to that of the Hessian (or inverse) on a Krylov basis of directions of
positive curvature. This choice guarantees positive definiteness of the estimate matrices.

Additionally, we present an implementation for these methods in Algorithm~\ref{alg:quNac} and a limited memory variant in Algorithm~\ref{alg:LquNac} in a Newton-CG framework. Both update variants exploit parallel linear algebra,
essentially performing multiple BFGS updates in parallel. This is apparently the first such parallel implementations of BFGS and L-BFGS updates.
Quadratic hereditary is proved for the full memory implementation.
Tests of linear SVM problems with Logistic Loss and a regularizer have shown the
\texttt{inverse} quNac method to be very promising, while our tests on Classic academic problems indicate that it is
robust. Certainly more exhaustive tests are required.

The flexibility afforded by the action constraint could potentially be used to incorporate 
these methods into various optimization frameworks, such as active set 
methods where the sampling matrix is the basis of kernel of active linear constraints.
Furthermore, using positive curvature is not the only possibility. Directions of negative curvature could be explored in a trust region
model~\cite{Gould2000,More1979a}.

% In particular, the flexibility 
%afforded by the action constraint could potentially be used to incorporate 
%these methods into various optimization frameworks, such as active set 
%methods where the sample matrix is the basis of kernel of active linear constraints, 
%or specialize our results to structured systems, such as KKT and saddle 
%point systems.

%or specialize our results to structured systems, such as KKT and saddle 
%point systems. 

%       The real merits of this family of estimate matrices is in its capacity to mesh with different optimization frameworks...
%Continuation methods such as Interior point where one carefully controls changes in a local quadratic model are another
%possibility, where the quNac methods could be used to track the changing inverse. Finally, in active set methods the
%action constraint could incorporate the null space of the current active set.

{\bf Acknowledgements and funding:}
The authors would like to thank Felix Lieder for his suggestions on constructing positive definite
matrices, and Artur Gower for proof reading the manuscript.

%%%%%%%%%%%%%%%%%%%%%%%%%%%%%%%%%%%%%%%%%%%%%
\parskip0pt
\baselineskip4pt
\bibliographystyle{siam}
% \printbibliography
\bibliography{ref/Automatic_Differentiation,ref/Optimization_methods,ref/quasi-Newton_Conjugate_Gradients,ref/Image_and_signal_recovery,ref/Graphs_Networks,ref/Preconditioners,ref/Now,ref/PDE_Constrained_Optimal_Control}

\begin{thebibliography}{10}

\bibitem{ginzburg-landau}
{\sc Igor~S Aranson and Lorenz Kramer}, {\em {The world of the complex
  Ginzburg-Landau equation}}, Reviews of Modern Physics, 74 (2002),
  pp.~99--143.

\bibitem{Bellavia2014}
{\sc Stefania Bellavia, Valentina~De Simone, Benedetta Morini, and Daniela
  di~Serafino}, {\em {On the update of constraint preconditioners for
  regularized KKT systems}}, Optimization Online,  (2014).

\bibitem{Inexact2006}
{\sc L.~Bergamaschi, R.~Bru, A.~Mart\'{\i}nez, and M.~Putti}, {\em
  {Quasi-Newton preconditioners for the inexact Newton method}}, Electronic
  Transactions on Numerical Analysis, 23 (2006), pp.~76--87.

\bibitem{Birgin2007}
{\sc E.~G. Birgin and J.~M. Mart\'{\i}nez}, {\em {Structured minimal-memory
  inexact quasi-Newton method and secant preconditioners for augmented
  Lagrangian optimization}}, Computational Optimization and Applications, 39
  (2007), pp.~1--16.

\bibitem{Broyden1965}
{\sc CG~Broyden}, {\em {A class of methods for solving nonlinear simultaneous
  equations}}, Mathematics of computation, 19 (1965), pp.~577--593.

\bibitem{Broyden1970}
{\sc C.~G. Broyden}, {\em {The Convergence of a Class of Double-rank
  Minimization Algorithms 1. General Considerations}}, J. Inst. Maths Applics,
  76 (1970), pp.~76--90.

\bibitem{Candes2009}
{\sc Emmanuel~J. Cand\`{e}s and Benjamin Recht}, {\em {Exact Matrix Completion
  via Convex Optimization}}, Foundations of Computational Mathematics, 9
  (2009), pp.~717--772.

\bibitem{Chang2011}
{\sc Chih-Chung Chang and Chih-Jen Lin}, {\em {Libsvm}}, ACM Transactions on
  Intelligent Systems and Technology, 2 (2011), pp.~1--27.

\bibitem{Christianson:1992}
{\sc Bruce Christianson}, {\em {Automatic Hessians by reverse accumulation}},
  IMA J. Numer. Anal., 12 (1992), pp.~135--150.

\bibitem{Davidon1959}
{\sc W.~C Davidon}, {\em {Variable metric method for minimization}}, tech.
  report, A.E.C. Research and Development Report, ANL-5990, 1959.

\bibitem{Dembo1982}
{\sc Ron~S. Dembo, Stanley~C. Eisenstat, and Trond Steihaug}, {\em {Inexact
  Newton Methods}}, SIAM Journal on Numerical Analysis, 19 (1982),
  pp.~400--408.

\bibitem{Dennis1979}
{\sc J.~E.~Jr. Dennis and R.~B. Schnabel}, {\em {Least Change Secant Updates
  for Quasi-Newton Methods}}, SIAM Review, 21 (1979), pp.~443--459.

\bibitem{Fletcher1960}
{\sc By~R Fletcher and M~J~D Powell}, {\em {A rapidly convergent descent method
  for minimization}}, The Computer Journal, 6 (1963), pp.~163----168.

\bibitem{Fletcher1970}
{\sc Rodger Fletcher}, {\em {A new approach to variable metric algorithms}},
  The Computer Journal, 13 (1970), pp.~317--323.

\bibitem{Fountoulakis2013}
{\sc Kimon Fountoulakis and Jacek Gondzio}, {\em {A Second-Order Method for
  Strongly Convex l1-regularization Problems}}, tech. report, Technical Report
  ERGO-13-011., 2013.

\bibitem{Gaul2012}
{\sc Andr\'{e} Gaul and Nico Schl\"{o}mer}, {\em {Preconditioned Recycling
  Krylov subspace methods for self-adjoint problems}}, ArXiv e-prints,  (2012),
  pp.~1--28.

\bibitem{Gill2001}
{\sc Philip~E Gill and Michael~W Leonard}, {\em {Reduced-Hessian quasi-Newton
  methods for unconstrained optimization}}, SIAM J. Optim., 12 (2001),
  pp.~209--237.

\bibitem{Giraud2000Inc}
{\sc L~Giraud, S~Gratton, and E~Martin}, {\em {Incremental spectral
  preconditioners for sequences of linear systems}}, Applied Numerical
  Mathematics, 57 (2007), pp.~1164--1180.

\bibitem{Goldfarb1970}
{\sc Donald Goldfarb}, {\em {A Family of Variable-Metric Methods Derived by
  Variational Means}}, Mathematics of Computation, 24 (1970), p.~23.

\bibitem{Gould2000}
{\sc N.~I.~M. Gould, S.~Lucidi, M.~Roma, and Ph. Toint}, {\em {Exploiting
  negative curvature directions in linesearch methods for unconstrained
  optimization}}, Optimization Methods and Software, 14 (2000), pp.~75--98.

\bibitem{Gower2014a}
{\sc R.~M. Gower}, {\em {Conjugate Gradients: The short and painful explanation
  with oblique projections}}, tech. report, University of Edinburgh, Maxwell
  Institute for Mathematical Sciences, 2014.

\bibitem{Gratton2011}
{\sc S~Gratton, A~Sartenaer, and J~Tshimanga}, {\em {On a class of limited
  memory preconditioners for large scale linear systems with multiple
  right-hand sides}}, SIAM Journal on Optimization, 21 (2011), pp.~912--935.

\bibitem{Greenstadt1969}
{\sc By~J Greenstadt}, {\em {Variations on Variable-Metric Methods}},
  Mathematics of Computation, 24 (1969), pp.~1--22.

\bibitem{Hestenes1952}
{\sc M.~R. Hestenes and E.~Stiefel}, {\em {Methods of Conjugate Gradients for
  Solving Linear Systems}}, Journal of research of the National Bureau of
  Standards, 49 (1952).

\bibitem{Huckle2007}
{\sc Thomas Huckle and Alexander Kallischko}, {\em {Frobenius Norm Minimization
  and Probing for Preconditioning}}, International Journal of Computer
  Mathematics, 00 (2007), pp.~1--31.

\bibitem{Loghin2006}
{\sc D.~Loghin, D.~Ruiz, and A.~Touhami}, {\em {Adaptive preconditioners for
  nonlinear systems of equations}}, Journal of Computational and Applied
  Mathematics, 189 (2006), pp.~362--374.

\bibitem{Mandel1993}
{\sc J~Mandel}, {\em {Balancing domain decomposition}}, Communications on
  Numerical Methods in Engineering, 9 (1993), pp.~233--241.

\bibitem{Morales2000a}
{\sc Jos\'{e}~Luis Morales and Jorge Nocedal}, {\em {Automatic Preconditioning
  by Limited Memory Quasi-Newton Updating}}, SIAM Journal on Optimization, 10
  (2000), pp.~1079--1096.

\bibitem{Mori1981}
{\sc JJ~Mor\'{e}, BS~Garbow, and KE~Hillstrom}, {\em {Testing unconstrained
  optimization software}}, ACM Transactions on Mathematical \ldots, 7 (1981),
  pp.~17--41.

\bibitem{More1979a}
{\sc JJ~Mor\'{e} and DC~Sorensen}, {\em {On the use of directions of negative
  curvature in a modified Newton method}}, Mathematical Programming, 16 (1979),
  pp.~1--20.

\bibitem{Nocedal1980}
{\sc Jorge Nocedal}, {\em {Updating Quasi-Newton Matrices with Limited
  Storage}}, Mathematics of Computation, 35 (1980), p.~773.

\bibitem{Nocedal1999}
{\sc J~Nocedal and S~J Wright}, {\em {Numerical Optimization}}, vol.~43 of
  Springer Series in Operations Research, Springer, 1999.

\bibitem{Parks2006}
{\sc Michael~L. Parks, Eric de~Sturler, Greg Mackey, Duane~D. Johnson, and
  Spandan Maiti}, {\em {Recycling Krylov Subspaces for Sequences of Linear
  Systems}}, SIAM Journal on Scientific Computing, 28 (2006), pp.~1651--1674.

\bibitem{Pearlmutter1993a}
{\sc Barak~A Pearlmutter}, {\em {Fast exact multiplication by the Hessian}},
  Tech. Report January, CSETech. Paper 286., 1993.

\bibitem{Petersen2012}
{\sc K.~B. Petersen and M.~S. Pedersen}, {\em {The Matrix Cookbook}}, tech.
  report, Technical University of Denmark, 2012.

\bibitem{Schnabel1983}
{\sc RB~Schnabel}, {\em {Quasi-Newton Methods Using Multiple Secant Equations;
  CU-CS-247-83}}, tech. report, Computer Science Technical Reports Computer.
  University of Colorado, Boulder, Boulder, 1983.

\bibitem{Shanno1971}
{\sc D~F Shanno}, {\em {Conditioning of Quasi-Newton Methods for Function
  Minimization}}, Mathematics of Computation, 24 (1971), pp.~647--656.

\bibitem{Shewchuk1994}
{\sc Jonathan~Richard Shewchuk}, {\em {An Introduction to the Conjugate
  Gradient Method Without the Agonizing Pain}}, tech. report, School of
  Computer Science Carnegie Mellon University, 1994.

\bibitem{sulem1999nonlinear}
{\sc C~Sulem and P~L Sulem}, {\em {The Nonlinear Schr\"{o}dinger Equation:
  Self-Focusing and Wave Collapse}}, no.~v. 139 in Applied Mathematical
  Sciences, Springer, 1999.

\bibitem{Tebbens2007}
{\sc Jurjen~Duintjer Tebbens}, {\em {Efficient Preconditioning of sequences of
  nonsymmetric linear systems}}, SIAM J. Matrix Anal. Appl., 29 (2007),
  pp.~1918--1941.

\bibitem{Woodbury1950}
{\sc Max~A Woodbury}, {\em {Inverting modified matrices}}, tech. report, Rep.
  no. 42, Statistical Research Group, Princeton University, 1950.

\bibitem{Yuan2012}
{\sc Guo-Xun Yuan, Chia-Hua Ho, and Chih-Jen Lin}, {\em {Recent Advances of
  Large-Scale Linear Classification}}, Proceedings of the IEEE, 100 (2012),
  pp.~2584--2603.

\end{thebibliography}

\section{Appendix: Updating the Inverse with the Direct approach} \label{sec:DFPinv}
Dispensing the iteration subscript $k$, to find the inverse $(G+E)^{-1}$ when a \texttt{direct} quNac update \ref{eq:genquNac}$(G,\mathcal{S}~\rightarrow~Q~\mathcal{S}).$
 is applied to $G$, we use the Woodbury formula~\cite{Woodbury1950}
\[(G+E)^{-1} = G^{-1}- G^{-1}U(I +VG^{-1}U)^{-1}V G^{-1},\]
where $G,E \in \R^{n\times n}$ and $E = UV$ with $U,V^T \in \R^{n \times q}.$
First we express the \texttt{direct} quNac update as  two rank-$p$ updates $G+E^1+E^2$
where
\begin{align}\label{eq:directquNacBapp}
(G+E) &= G+ \underbrace{(Q -G) \proj{\mathcal{S}}{ Q} Q}_{E^{1}} -\underbrace{Q\proj{\mathcal{S}}{ Q}G\left(I
-\proj{\mathcal{S}}{
Q} Q\right)}_{E^2},
\end{align}
The first $E^1$
 can be split up as $E^1 = U^1V^1$ with
\[U^1 = (Q -G) D, \quad V^1 = (\mathcal{S}^T Q \mathcal{S})^{-1}\mathcal{S}^T Q. \]
Applying the Woodbury formula where $H \equiv G^{-1}$ we get
\begin{align*}
(G+ E^1)^{-1} &= H- H (Q -G) \mathcal{S}\left(\phantom{\sum} \hspace{-0.5cm}I + (\mathcal{S}^T Q
\mathcal{S})^{-1}\mathcal{S}^T Q H (Q -G)
\mathcal{S}\right)^{-1}(\mathcal{S}^T Q \mathcal{S})^{-1}\mathcal{S}^T Q H\\
&= H- H (Q -G) \mathcal{S}\left(\phantom{\sum} \hspace{-0.5cm} (\mathcal{S}^T Q \mathcal{S})^{-1}\mathcal{S}^T Q H
Q \mathcal{S}\right)^{-1}(\mathcal{S}^T Q
\mathcal{S})^{-1}\mathcal{S}^T Q H\\
&=H- H (Q -G) \mathcal{S}\left(\mathcal{S}^T Q H Q \mathcal{S}\right)^{-1}\mathcal{S}^T Q H\\
&= H- H (Q -G) \proj{\mathcal{S}}{Q H Q}Q H.
\end{align*}
The second update can be split up as $E^2 = U^2 V^2$ with
\[U^2 = -Q \mathcal{S}(\mathcal{S}^T Q \mathcal{S})^{-1}= (V^1)^T , \quad V^2 = \mathcal{S}^T G\left(I
-\proj{\mathcal{S}}{
Q} Q\right). \]
If we let $\bar{H} = (G +E^1)^{-1}$, then applying the Woodbury formula again
\begin{align*}
&((G+E^1)+ E^2)^{-1} = \bar{H}\\
&+ \underbrace{\bar{H}Q \mathcal{S}(\mathcal{S}^T Q \mathcal{S})^{-1}}_{\mathbb{I}}\left(\underbrace{I
-\mathcal{S}^T
G\left(I -\proj{\mathcal{S}}{Q} Q\right)\bar{H}Q \mathcal{S}(\mathcal{S}^T
Q \mathcal{S}}_{\mathbb{II}})^{-1}\right)^{-1}\underbrace{\mathcal{S}^T G\left(I -\proj{\mathcal{S}}{
Q} Q\right) \bar{H}}_{\mathbb{III}}.
\end{align*}

When substituting in $\bar{H}$, simplifications arise such as
\begin{align*}\bar{H}Q \mathcal{S} &= \left( H- H (Q -G) \proj{\mathcal{S}}{Q H Q}Q H
\right)Q \mathcal{S} \\
&= \left( H Q \mathcal{S} - H (Q -G)\mathcal{S}\right) \\
&= \mathcal{S}.
\end{align*}

Thus
\begin{align*}
\mathbb{I}& =\bar{H}Q \mathcal{S}(\mathcal{S}^T Q \mathcal{S})^{-1} = \mathcal{S} (\mathcal{S}^T Q
\mathcal{S})^{-1},
\end{align*}
and
\begin{align*} \mathbb{II}& =I -\mathcal{S}^T G\left(I -\proj{\mathcal{S}}{
Q} Q\right)\bar{H}Q \mathcal{S}(\mathcal{S}^T Q \mathcal{S})^{-1} \\
& =I -\mathcal{S}^T G\left(I -\proj{\mathcal{S}}{
Q} Q\right)D(\mathcal{S}^T Q \mathcal{S})^{-1}\\
&= I -\mathcal{S}^TG(\mathcal{S} -\mathcal{S}) (\mathcal{S}^T Q \mathcal{S})^{-1} = I.
%\mathcal{S}^T G\proj{\mathcal{S}}{
%Q} Q  D(\mathcal{S}^T Q \mathcal{S})^{-1} =I.
\end{align*}
For the final part, take note that

\begin{align*}\mathcal{S}^T Q\bar{H} &= \mathcal{S}^T Q\left( H- H (Q -G) \proj{\mathcal{S}}{Q H
Q}Q H \right) \\
&= \mathcal{S}^T Q H-  \mathcal{S}^T Q(HQ -I) \proj{\mathcal{S}}{Q H Q}Q H  \\
&= \mathcal{S}^T Q H + \mathcal{S}^T Q \proj{\mathcal{S}}{Q H Q}Q H -\mathcal{S}^T Q H  \\
&= \mathcal{S}^T Q \proj{\mathcal{S}}{Q H Q}Q H.
\end{align*}
Furthermore
\begin{align*}\mathcal{S}^T G\bar{H} &= \mathcal{S}^T G\left( H- H (Q -G) \proj{\mathcal{S}}{Q H Q}Q H
\right) \\
&= \mathcal{S}^T \left(I +   (G-Q) \proj{\mathcal{S}}{Q H Q}Q H \right)
\end{align*}

Thus

\begin{align*}
\mathbb{III}& =\mathcal{S}^T G\left(I -\proj{\mathcal{S}}{Q} Q\right) \bar{H} \\
&= \mathcal{S}^T G \bar{H}  -\mathcal{S}^T G D (\mathcal{S}^TQ \mathcal{S} )^{-1} \mathcal{S}^T Q \bar{H} \\
&=\mathcal{S}^T G \bar{H}  - \mathcal{S}^T G \proj{\mathcal{S}}{Q H
Q}Q H \\
&=  \mathcal{S}^T \left(I +   (G-Q) \proj{\mathcal{S}}{Q H Q}Q H \right)-\mathcal{S}^T G
\proj{\mathcal{S}}{Q H Q}  Q H\\
&= \mathcal{S}^T\left(I    -Q\proj{\mathcal{S}}{Q H Q}Q H \right).
\end{align*}

Bringing all this together yields

\begin{align*} (G+E)^{-1} &=
\overbrace{H- H (Q -G) \proj{\mathcal{S}}{Q H Q}Q H}^{\bar{H}}+ \overbrace{\proj{\mathcal{S}}{Q}\left(I
-Q\proj{\mathcal{S}}{Q H Q}Q H \right)}^{\mathbb{I} \, \cdot \, \mathbb{II} \, \cdot \, \mathbb{III}} \nonumber
\\
&=  H-  (HQ -I) \proj{\mathcal{S}}{Q H Q}Q H+ \proj{\mathcal{S}}{Q}-
\proj{\mathcal{S}}{Q H Q}Q H \nonumber \\
&=  H+ \proj{\mathcal{S}}{Q}-  HQ \proj{\mathcal{S}}{Q H Q}Q H. \label{eq:DFPinv}
\end{align*}

With indices
\begin{equation}
\label{eq:DFPinv}(G_k+E_k)^{-1} =H_k+ \proj{\mathcal{S}_k}{Q_{k+1}}-  H_kQ_{k+1} \proj{\mathcal{S}_k}{Q_{k+1} H_k
Q_{k+1}}Q_{k+1}
H_k.
\end{equation}

\end{document}